\documentclass[dvips,preprint]{imsart}

\usepackage{amsmath}
\usepackage{amsfonts}
\usepackage{mathrsfs}
\usepackage{color}
\usepackage{ amsmath, amsfonts, amssymb, amsthm, amscd}
\usepackage[T1]{fontenc}
\usepackage[english]{babel}

\DeclareSymbolFont{calletters}{OMS}{cmsy}{m}{n}
\DeclareSymbolFontAlphabet{\mathcal}{calletters}

\def\be{\begin{eqnarray}}
\def\ee{\end{eqnarray}}

\def\b*{\begin{eqnarray*}}
\def\e*{\end{eqnarray*}}

\newtheorem{Theorem}{Theorem}[section]
\newtheorem{Definition}[Theorem]{Definition}
\newtheorem{Proposition}[Theorem]{Proposition}

\newtheorem{Assumption}[Theorem]{Assumption}
\newtheorem{Lemma}[Theorem]{Lemma}

\newtheorem{Remark}[Theorem]{Remark}

\newcommand{\No}[1]{\left\|#1\right\|}     
\newcommand{\abs}[1]{\left|#1\right|}     

\def \D{\mathbb{D}}
\def \E{\mathbb{E}}
\def \F{\mathbb{F}}
\def \H{\mathbb{H}}

\def \P{\mathbb{P}}

\def\Pc{{\cal P}}

\def\Fh{{\hat F}}

\def\Tr#1{{\rm Tr}\left[#1\right]}

\def \Sup{\displaystyle\sup}

\def\einf{{\rm ess \, inf}}
\def\esup{{\rm ess \, sup}}
\def\trace{{\rm Tr}}

\def\={\;=\;}
\def\.{\;.}

\def\eps{\varepsilon}

\def\reff#1{{\rm(\ref{#1})}}

\def\1{{\bf 1}}
\def \ep{\hbox{ }\hfill{ ${\cal t}$~\hspace{-5.1mm}~${\cal u}$   } }

\def \proof{{\noindent \bf Proof. }}
\def \ep{\hbox{ }\hfill$\Box$}

 \def\normeL2#1{\left\|{#1}\right\|_{L^2}}

\begin{document}
\begin{frontmatter}
 \title{Second-order BSDEs with general reflection and game options under uncertainty}

\date{}
\runtitle{}
\author{\fnms{Anis}
 \snm{MATOUSSI}\corref{}\ead[label=e1]{anis.matoussi@univ-lemans.fr}}
\address{
 LUNAM Universit\'e, Universit\'e du Maine\\
 F\'ed\'eration de Recherche 2962 du CNRS\\
 Math\'ematiques des Pays de Loire \\
Laboratoire Manceau de Math\'ematiques\\
Avenue Olivier Messiaen\\ F-72085 Le Mans Cedex 9, France\\
 and
CMAP,  Ecole Polytechnique, Palaiseau, France\\\printead{e1}}

\author{\fnms{Lambert}
 \snm{PIOZIN}\corref{}\ead[label=e2]{lambert.piozin.etu@univ-lemans.fr}}
\address{
Laboratoire Manceau de Math\'ematiques\\
 Avenue Olivier Messiaen\\ F-72085 Le Mans Cedex 9, France\\\printead{e2}}

\author{\fnms{Dylan}
 \snm{POSSAMA\"{I}}\corref{}\ead[label=e3]{possamai@ceremade.dauphine.fr}}
\thankstext{t2}{Research partly supported by the Chair {\it Financial Risks} of the {\it Risk Foundation} sponsored by Soci\'et\'e G\'en\'erale, the Chair {\it Derivatives of the Future} sponsored by the {F\'ed\'eration Bancaire Fran\c{c}aise}, and the Chair {\it Finance and Sustainable Development} sponsored by EDF and Calyon }
\address{
Universit\'e Paris-Dauphine\\
 Ceremade, bureau B518,\\
  Place du Mar\'echal de Lattre de Tassigny,\\
   75775 Paris Cedex 16, France \\\printead{e3}}

\runauthor{A. Matoussi,  L. Piozin,  D. Possama\"{i}}

\vspace{3mm}

 \begin{abstract} : The aim of this paper is twofold. First, we extend the results of \cite{mpz} concerning the existence and uniqueness of second-order reflected 2BSDEs to the case of two obstacles. Under some regularity assumptions on one of the barriers, similar to the ones in \cite{crep}, and when the two barriers are completely separated, we provide a complete wellposedness theory for doubly reflected second-order BSDEs. We also show that these objects are related to non-standard optimal stopping games, thus generalizing the connection between DRBSDEs and Dynkin games first proved by Cvitani\'c and Karatzas \cite{cvi}. More precisely, we show under a technical assumption that the second order DRBSDEs provide solutions of what we call uncertain Dynkin games and that they also allow us to obtain super and subhedging prices for American game options (also called Israeli options) in financial markets with volatility uncertainty.\footnote{The authors are grateful to Chao Zhou, Jianfeng Zhang and Marcel Nutz for their interesting comments and remarks . They also would like to thank warmly two anonymous referees and an associate editor who greatly helped to improve this paper.}

\vspace{10mm}

\vspace{5mm}

\end{abstract}
\end{frontmatter}

\section{Introduction}
\label{intro}
The theory of Backward Stochastic Differential Equations (BSDEs for short) was initiated by Bismut in \cite{bis}, where their linear version appeared as an equation for the adjoint process in the stochastic version of Pontryagin maximum principle. This notion was then generalized by Pardoux and Peng in \cite{pardpeng}, where they were the first to consider these equations in their full generality. Since then BSDEs have been widely used in stochastic control and especially in mathematical
finance, as any pricing problem by replication can be written in terms of linear BSDEs, or non-linear
BSDEs when portfolios constraints are taken into account as in El Karoui, Peng and Quenez \cite{elkaroui}.

\vspace{0.35em}
On a filtered probability space $(\Omega, \mathbb F,(\mathcal F_t)_{0\leq t\leq T},\mathbb P)$  carrying a Brownian motion $B$, solving a BSDE consists in finding a couple $(Y,Z)$ of progressively measurable processes such that
$$Y_t=\xi +\int_t^Tf_s(Y_s,Z_s)ds-\int_t^TZ_sdB_s,\text{ } t\in [0,T], \text{ }\mathbb P-a.s.,$$

where $f$ (also called the generator) is progressively measurable and $\xi$ (the terminal condition) is a $\mathcal{F}_T$-measurable random variable.  When $f$ satisfies a uniform Lipschitz assumption and both $f$ and $\xi$ are square-integrable, Pardoux and Peng \cite{pardpeng} obtained a wellposedness theory for these equations.

\vspace{0.35em}
A few years later, the five authors El Karoui, Kapoudjian,  Pardoux, Peng and  Quenez in
\cite{elkarkap}, introduced the notion of reflected BSDEs (RBSDEs for short). Now, in addition to the generator $f$ and the terminal condition $\xi$, they also consider another progressively measurable process $S$ which will play the role of a barrier. In that framework, we look for a triple of progressively measurable processes $(Y,Z,K)$, where $K$ is furthermore non-decreasing, such that
\begin{align*}
&Y_t=\xi+ \int_t^T f_s(Y_s,Z_s)ds+K_T-K_t-\int_t^T Z_sdB_s, \,\,t\in[0,T],\ \mathbb P-a.s.\\
&Y_t\geq S_t , \ t\in[0,T],\ \mathbb P-a.s.\\
&\int_0^T(Y_t-S_t)dK_t=0,\ \mathbb P-a.s.
\end{align*}
The role of the process $K$ here is to push upward $Y$ in
order to keep it above the barrier $S$. The last condition is known as the Skorohod condition and guarantees that the process $K$ acts in a minimal way, that is to say only when the process $Y$ reaches the lower barrier $S$. The development of RBSDEs has been motivated
by the problem of pricing American contingent claim by replication, especially in
constrained markets (see \cite{elkarquen} and \cite{elkarmat} for more details), and for their natural connection with variational inequalities and the obstacle problem for deterministic and/or stochastic quasilinear PDEs (see  \cite{bally}, \cite{elkarkap} and \cite{MS10}).

\vspace{0.35em}
Building upon these results, Cvitani\'c and Karatzas \cite{cvi} have then introduced the notion of BSDEs
with two reflecting barriers. Roughly speaking, in \cite{cvi} (see also \cite{elkarmat} and \cite{HLM} to name but a few) the authors have looked for a solution to a BSDE whose $Y$ component is forced to stay between two prescribed processes $L$ and $
U$, $(L\leq U)$. More precisely, they were looking for a quadruple of progressively measurable processes $(Y,Z,K^+,K^-)$, where $K^+$ and $K^-$ are in addition non-decreasing such that
\begin{align*}
&Y_t=\xi+ \int_t^T f_s(Y_s,Z_s)ds+K^-_T-K^-_t-K_T^++K_t^+-\int_t^T Z_sdB_s, \,\,t\in[0,T],\ \mathbb P-a.s.\\
&L_t\leq Y_t\leq U_t , \ t\in[0,T],\ \mathbb P-a.s.\\
&\int_0^T(Y_t-L_t)dK^-_t=\int_0^T(U_t-Y_t)dK^+_t=0,\ \mathbb P-a.s.
\end{align*}

These BSDEs have been developed especially in connection
with Dynkin games, mixed differential games and {\it{recallable}}
options (see (\cite{elkarmat}, \cite{Hamadene95}, \cite{Lepeltier95} and \cite{ham2}). It is now established that under quite general assumptions, including in models with
jumps, existence of a solution to a (simply) reflected BSDE is guaranteed under mild conditions, whereas existence of a solution to a \emph{doubly reflected} BSDE (DRBSDE for short) is equivalent to the so-called {\em Mokobodski condition}. This condition essentially postulates the existence of a quasimartingale between the barriers (see in particular \cite{HH06}). As for uniqueness of solutions, it is guaranteed under mild integrability conditions (see e.g. \cite[Remark 4.1]{HH06}). However, for practical purposes, existence and uniqueness are not the only relevant results and may not be enough. For instance, one can consider the problem of pricing convertible bonds in finance using the DRBSDE theory (see \cite{BCJR1}, \cite{BCJR2} and \cite{crep}). In this case, the state-process (first component) $Y$ may be interpreted in terms of an arbitrage price process for the bond. As demonstrated in \cite{BCJR2}, the mere existence of a solution to the related DRBSDE is a result with important theoretical consequences in terms of pricing and hedging the bond. Yet, in order to give further developments to these results in Markovian set-ups, Cr\'epey and Matoussi \cite{crep} have established \emph{bound and error estimates} and {\em comparison theorem} for DRBSDE, which require more regularity assumptions on the barriers.

\vspace{0.35em}
More recently, motivated by numerical methods for fully nonlinear PDEs, second order BSDEs (2BSDEs for short) were introduced by Cheredito, Soner, Touzi and Victoir in \cite{cstv}. Then Soner, Touzi, Zhang \cite{stz} proposed a new formulation and obtained a complete theory of existence and uniqueness for such BSDEs. The main novelty in their approach is that they require that the solution verifies the equation $\mathbb P-a.s.$ for every probability measure $\mathbb P$ in a non-dominated set. Their approach therefore shares many connections with the deep theory of quasi-sure analysis initiated by Denis and Martini \cite{denis} and the $G$-expectations developed by Peng \cite{peng}.

\vspace{0.35em}
Intuitively speaking (we refer the reader to \cite{stz} for more details), the solution to a 2BSDE with generator $F$ and terminal condition $\xi$ can be understood as a supremum in some sense of the classical BSDEs with the same generator and terminal condition, but written under the different probability measures considered. Following this intuition, a non-decreasing process $K$ is added to the solution and it somehow pushes (in a minimal way) the solution so that it stays above the solutions of the classical BSDEs.

\vspace{0.35em}
Following these results and motivated by the pricing of American contingent claims in markets with volatility uncertainty, Matoussi, Possama\"i and Zhou \cite{mpz} used the methodology of \cite{stz} to introduce a notion of reflected second order BSDEs, and proved existence and uniqueness in the case of a lower obstacle. The fact that they consider only lower obstacles was absolutely crucial. Indeed, as mentioned above, in that case, the effects due to the reflection and the second order act in the same direction, in the sense that they both force the solution to stay above some processes. One therefore only needs to add a non-decreasing process to the solution of the equation. However, as soon as one tries to consider upper obstacles, the two effects start to counterbalance each other and the situation changes drastically. This case was thus left open in \cite{mpz}. On a related note, we would like to refer the reader to the very recent article \cite{etz}, which gives some specific results for the optimal stopping problem under a non-linear expectation (which roughly corresponds to a 2RBSDE with generator equal to $0$). However, since it is a "sup-sup" problem, it is only related to the lower reflected 2BSDEs. Even more recently and after the completion of this paper, Nutz and Zhang \cite{nutz} managed to treat the same problem of optimal stopping under non-linear expectations but now with an "inf-sup" formulation, which, as shown by Proposition \ref{dynkin}, is related to upper reflected 2BSDEs. 

\vspace{0.35em}
The first aim of this paper is to extend the results of \cite{mpz} to the case of doubly reflected second-order BSDEs when we assume enough regularity on one of the barriers (as in \cite{crep}) and that the two barriers are completely separated (as in \cite{HH05} and \cite{HH06}). In that case, we show that the right way to define a solution is to consider a 2BSDE where we add a process $V$ which has only bounded variations (see Definition \ref{death2}). Our next step towards a theory of existence and uniqueness is then to understand as much as possible how and when this bounded variation process acts. Our key result is obtained in Proposition \ref{prop.imp}, and allows us to obtain a special Jordan decomposition for $V$, in the sense that we can decompose it into the difference of two non-decreasing processes which never act at the same time. Thanks to this result, we are then able to obtain a priori estimates and a uniqueness result. Next, we reuse the methodology of \cite{mpz} to construct a solution.

\vspace{0.35em}
We also show that these objects are related to non-standard optimal stopping games, thus generalizing the connection between DRBSDEs and Dynkin games first proved by Cvitani\'c and Karatzas \cite{cvi}. Finally, we show that the second order DRBSDEs allow to obtain super and subhedging prices for American game options (also called Israeli options) in financial markets with volatility uncertainty and that, under a technical assumption, they provide solutions of what we call uncertain Dynkin games.

\vspace{0.35em}
The paper is organized as follows. After recalling some notations and definitions, in Section \ref{section.1}, we treat the problem of uniqueness in Section \ref{section.2}. Section \ref{section.3} is then devoted to the pathwise construction of a solution, thus solving the existence problem. Finally, we investigate in Section \ref{sec.app} the aforementioned game theoretical and financial applications. The Appendix is devoted to some technical results used throughout the paper.

\section{Definitions and Notations}\label{section.1}

\subsection{The stochastic framework}

Let $\Omega:=\left\{\omega\in C([0,T],\mathbb R^d):\omega_0=0\right\}$ be the canonical space equipped with the uniform norm $\No{\omega}_{\infty}:=\sup_{0\leq t\leq T}|\omega_t|$, $B$ the canonical process, $\mathbb P_0$ the Wiener measure, $\mathbb F:=\left\{\mathcal F_t\right\}_{0\leq t\leq T}$ the filtration generated by $B$, and $\mathbb F^+:=\left\{\mathcal F_t^+\right\}_{0\leq t\leq T}$ the right limit of $\mathbb F$. A probability measure $\mathbb P$ will be called a local martingale measure if the canonical process $B$ is a local martingale under $\mathbb P$. Then, using results of Bichteler \cite{bich} (see also Karandikar \cite{kar} for a modern account), the quadratic variation $\langle B\rangle$ and its density $\widehat a$ can be defined pathwise, and such that they coincide with the usual definitions under any local martingale measure. 

\vspace{0.35em}
With the intuition of modeling volatility uncertainty, we let $\overline{\mathcal P}_W$ denote the set of all local martingale measures $\mathbb P$ such that
\begin{equation}
\langle B\rangle \text{ is absolutely continuous in $t$ and $\widehat a$ takes values in $\mathbb S_d^{>0}$, }\mathbb P-a.s.,
\label{eq:}
\end{equation}
where $\mathbb S_d^{>0}$ denotes the space of all $d\times d$ real valued positive definite matrices. 

\vspace{0.5em}
However, since this set is too large for our purpose (in particular there are examples of measures in $\overline{\Pc}_W$ which do not satisfy the martingale representation property, see \cite{stz3} for more details), we will concentrate on the following subclass $\Pc_S$ consisting of
\begin{equation}
\mathbb P^\alpha:=\mathbb P_0\circ (X^\alpha)^{-1} \text{ where }X_t^\alpha:=\int_0^t\alpha_s^{1/2}dB_s,\text{ }t\in [0,T],\text{ }\mathbb P_0-a.s.,
\end{equation}
for some $\mathbb F$-progressively measurable process $\alpha$ taking values in $\mathbb S_d^{>0}$ with $\int_0^T|\alpha_t|dt<+\infty$, $\mathbb {P}_0-a.s.$

\vspace{0.35em}
This subset has the convenient property that all its elements do satisfy the martingale representation property and the Blumenthal $0-1$ law (see \cite{stz3} for details), which are crucial tools for the BSDE theory.

\subsection{Generator and measures}
We consider a map $H_t(\omega,y,z,\gamma):[0,T]\times\Omega\times\mathbb{R}\times\mathbb{R}^d\times D_H\rightarrow \mathbb{R}$, where $D_H \subset \mathbb{R}^{d\times d}$ is a given subset containing $0$, whose Fenchel transform w.r.t. $\gamma$ is denoted by
\begin{align*}
&F_t(\omega,y,z,a):=\underset{\gamma \in D_H}{\Sup}\left\{\frac12\trace(a\gamma)-H_t(\omega,y,z,\gamma)\right\} \text{ for } a \in \mathbb S_d^{>0},\\[0.3em]
&\widehat{F}_t(y,z):=F_t(y,z,\widehat{a}_t) \text{ and } \widehat{F}_t^0:=\widehat{F}_t(0,0).
\end{align*}

We denote by $D_{F_t(y,z)}:=\left\{a, \text{ }F_t(\omega,y,z,a)<+\infty\right\}$ the domain of $F$ in $a$ for a fixed $(t,\omega,y,z)$. As in \cite{stz} we fix a constant $\kappa \in (1,2]$ and restrict the probability measures in $\mathcal{P}_H^\kappa\subset \overline{\mathcal{P}}_S$

\begin{Definition}\label{def}
$\mathcal{P}_H^\kappa$ consists of all $\mathbb P \in \overline{\mathcal{P}}_S$ such that
$$\underline{a}^\mathbb P\leq \widehat{a}\leq \bar{a}^\mathbb P, \text{ } dt\times d\mathbb P-a.s. \text{ for some } \underline{a}^\P, \bar{a}^\P \in \mathbb S_d^{>0}, \text{ and } \phi^{2,\kappa}_H<+\infty,
$$
where
$$\phi^{2,\kappa}_H:=\underset{\mathbb P\in\mathcal P^\kappa_H}{\sup}\mathbb E^{\mathbb P}\left[\underset{0\leq t\leq T}{\esup}^{\mathbb P}\left(\mathbb E_t^{H,\mathbb P}\left[\int^T_0|\Fh^0_s|^\kappa ds\right]\right)^{\frac{2}{\kappa}}\right].$$
\end{Definition}

\vspace{0.35em}
\begin{Definition}
We say that a property holds $\mathcal{P}^\kappa_H$-quasi-surely ($\mathcal{P}^\kappa_H$-q.s. for short) if it holds $\mathbb P $-a.s. for all $\mathbb P\in \mathcal{P}^\kappa_H$.
\end{Definition}

We now state the main assumptions on the function $F$ which will be our main interest in the sequel
\begin{Assumption} \label{assump.href}
\begin{itemize}
\item[\rm{(i)}] The domain $D_{F_t(y,z)}=D_{F_t}$ is independent of $(\omega,y,z)$.

\vspace{0.35em}
\item[\rm{(ii)}] For fixed $(y,z,a)$, $F$ is $\mathbb{F}$-progressively measurable in $D_{F_t}$.

\vspace{0.35em}
\item[\rm{(iii)}] We have the following uniform Lipschitz-type property in $y$ and $z$
$$\forall (y,y',z,z',t,a,\omega), \text{ } \abs{ F_t(\omega,y,z,a)- F_t(\omega,y',z',a)}\leq C\left(\abs{y-y'}+\abs{ a^{1/2}\left(z-z'\right)}\right).$$
\item[\rm{(iv)}] $F$ is uniformly continuous in $\omega$ for the $||\cdot||_\infty$ norm.

\vspace{0.35em}
\item[\rm{(v)}] $\Pc^\kappa_H$ is not empty. 
\end{itemize}
\end{Assumption}

\begin{Remark}
Assumptions (ii), (iii) are completely standard in the BSDE literature since the paper \cite{pardpeng}. Similarly, (i) was already present in the first paper on 2BSDEs in a quasi-sure formulation \cite{stz} and is linked to the fact that one does not know how to treat coupled second-order FBSDEs. The last hypothesis (iv) is also proper to the second order framework, and allows us to not only give a pathwise construction for the solution to the 2RBSDE, but to recover the very important dynamic programming property. We refer the reader to Section \ref{section.3} for more details.
\end{Remark}

%

\subsection{Quasi-sure norms and spaces}

The following spaces and the corresponding norms will be used throughout the paper. With the exception of the space $\mathbb L^{p,\kappa}_H$, they all are immediate extensions of the usual spaces to the quasi-sure setting.

\vspace{0.35em}
For $p\geq 1$, $L^{p,\kappa}_H$ denotes the space of all $\mathcal F_T$-measurable scalar r.v. $\xi$ with
$$\No{\xi}_{L^{p,\kappa}_H}^p:=\underset{\mathbb{P} \in \mathcal{P}_H^\kappa}{\sup}\mathbb E^{\mathbb P}\left[|\xi|^p\right]<+\infty.$$

$\mathbb H^{p,\kappa}_H$ denotes the space of all $\mathbb F^+$-progressively measurable $\mathbb R^d$-valued processes $Z$ with
$$\No{Z}_{\mathbb H^{p,\kappa}_H}^p:=\underset{\mathbb{P} \in \mathcal{P}_H^\kappa}{\sup}\mathbb E^{\mathbb P}\left[\left(\int_0^T|\widehat a_t^{1/2}Z_t|^2dt\right)^{\frac p2}\right]<+\infty.$$

$\mathbb D^{p,\kappa}_H$ denotes the space of all $\mathbb F^+$-progressively measurable $\mathbb R$-valued processes $Y$ with
$$\mathcal P^\kappa_H-q.s. \text{ c\`adl\`ag paths, and }\No{Y}_{\mathbb D^{p,\kappa}_H}^p:=\underset{\mathbb{P} \in \mathcal{P}_H^\kappa}{\sup}\mathbb E^{\mathbb P}\left[\underset{0\leq t\leq T}{\sup}|Y_t|^p\right]<+\infty,$$
where c\`adl\`ag is the french acronym for "right-continuous with left-limits".

\vspace{0.35em}
$\mathbb I^{p,\kappa}_H$ denotes the space of all $\mathbb F^+$-progressively measurable $\mathbb R$-valued processes $K$ null at $0$ with
$$\mathcal P^\kappa_H-q.s., \text{ c\`adl\`ag and non-decreasing paths, and }\No{K}_{\mathbb I^{p,\kappa}_H}^p:=\underset{\mathbb{P} \in \mathcal{P}_H^\kappa}{\sup}\mathbb E^{\mathbb P}\left[\left(K_T\right)^p\right]<+\infty.$$

$\mathbb V^{p,\kappa}_H$ denotes the space of all $\mathbb F^+$-progressively measurable $\mathbb R$-valued processes $V$ null at $0$ with paths which are $\mathcal P^\kappa_H-q.s. \text{ c\`adl\`ag and of bounded variation, and such that}$
$$\No{V}_{\mathbb V^{p,\kappa}_H}^p:=\underset{\mathbb{P} \in \mathcal{P}_H^\kappa}{\sup}\mathbb E^{\mathbb P}\left[\left( \rm{Var}_{0,T}  (V)\right)^p\right]<+\infty.$$

For each $\xi \in L^{1,\kappa}_H$, $\mathbb P\in \mathcal P^\kappa_H$ and $t \in [0,T]$ denote
$$\mathbb E_t^{H,\mathbb P}[\xi]:=\underset{\mathbb P^{'}\in \mathcal P^\kappa_H(t^{+},\mathbb P)}{\esup^{\mathbb P}}\mathbb E^{\mathbb P^{'}}_t[\xi] \text{ where } \mathcal P^\kappa_H(t^{+},\mathbb P):=\left\{\mathbb P^{'}\in\mathcal P^\kappa_H:\mathbb P^{'}=\mathbb P \text{ on }\mathcal F_t^+\right\}.$$

Here $\mathbb{E}_t^{\mathbb P}[\xi]:=E^{\mathbb P}[\xi|\mathcal F_t]$. Then we define for each $p\geq \kappa$,
$$\mathbb L_H^{p,\kappa}:=\left\{\xi\in L^{p,\kappa}_H:\No{\xi}_{\mathbb L_H^{p,\kappa}}<+\infty\right\} \text{ where } \No{\xi}_{\mathbb L_H^{p,\kappa}}^p:=\underset{\mathbb P\in\mathcal P^\kappa_H}{\sup}\mathbb E^{\mathbb P}\left[\underset{0\leq t\leq T}{\esup}^{\mathbb P}\left(\mathbb E_t^{H,\mathbb P}[|\xi|^\kappa]\right)^{\frac{p}{\kappa}}\right].$$

We denote by $\mbox{UC}_b(\Omega)$ the collection of all bounded and uniformly continuous maps $\xi:\Omega\rightarrow \mathbb R$ with respect to the $\No{\cdot}_{\infty}$-norm, and we let
$$\mathcal L^{p,\kappa}_H:=\text{the closure of $\mbox{UC}_b(\Omega)$ under the norm $\No{\cdot}_{\mathbb L^{p,\kappa}_H}$, for every $1\leq\kappa \leq p$}.$$

\vspace{0.35em}
Finally, for every $\P\in\Pc^\kappa_H$, and for any $p\geq 1$, $L^p(\P)$, $\H^p(\P)$, $\D^p(\P)$, $\mathbb I^p(\P)$ and $\mathbb V^p(\P)$ will denote the corresponding usual spaces when there is only one measure $\P$.

\subsection{Obstacles and definition}

First, we consider a process $S$ which will play the role of the upper obstacle. We will always assume that $S$ verifies the following properties

\vspace{0.35em}
\begin{Assumption}\label{assump.obs}
$\rm{(i)}$ $S$ is $\mathbb F$-progressively measurable and c\`adl\`ag.

\vspace{0.35em}
$\rm{(ii)}$ $S$ is uniformly continuous in $\omega$ in the sense that for all $t$
$$\abs{S_t(\omega)-S_t(\widetilde{\omega})}\leq\rho\left(\No{\omega-\widetilde{\omega}}_t\right),\text{ }\forall \text{ }(\omega,\widetilde{\omega})\in\Omega^2,$$
for some modulus of continuity $\rho$ and where we define $\No{\omega}_t:=\underset{0\leq s\leq t}{\sup}\abs{\omega(s)}$.

\vspace{0.35em}
$\rm{(iii)}$ $S$ is a semimartingale for every $\mathbb P\in\mathcal P^\kappa_H$, with the decomposition
\begin{equation}
S_t=S_0+\int_0^tP_sdB_s+A_t^\mathbb P,\ \mathbb P-a.s.,\text{ for all } \mathbb P\in\mathcal P^\kappa_H,
\label{eq:decomp}
\end{equation}
where the $A^\mathbb P$ are bounded variation processes with Jordan decomposition $A^{\mathbb P,+}-A^{\mathbb P,-}$ and
$$\zeta^{2,\kappa}_H:=\underset{\mathbb P\in\mathcal P^\kappa_H}{\sup}\left(\mathbb E^\mathbb P\left[\underset{0\leq t\leq T}{\esup^\mathbb P}\ \mathbb E_t^{H,\mathbb P}\left[\left(\int_t^T\abs{\widehat a_s^{1/2}P_s}^2ds\right)^{\kappa/2}+\left(A_T^{\mathbb P,+}\right)^\kappa\right]\right]\right)^2<+\infty.$$

\vspace{0.35em}
$\rm{(iv)}$ $S$ satisfies the following integrability condition
\begin{align*}
\psi^{2,\kappa}_H&:=\underset{\mathbb P\in\mathcal P^\kappa_H}{\sup}\mathbb E^{\mathbb P}\left[\underset{0\leq t\leq T}{\esup}^{\mathbb P}\left(\mathbb E_t^{H,\mathbb P}\left[\underset{0\leq s\leq T}{\sup}\abs{S_s}^{\kappa}\right]\right)^{\frac{2}{\kappa}}\right]<+\infty.
\end{align*}
\end{Assumption}

\vspace{0.35em}
\begin{Remark}
We assumed here that $S$ was a semimartingale. This is directly linked to the fact that this is one of the conditions under which existence and uniqueness of solutions to standard doubly reflected BSDEs with upper obstacle $S$ are guaranteed. More precisely, this assumption is needed for us in the proof of Lemma \reff{lemd}, and it will be also crucial in order to obtain a priori estimates for 2BSDEs with two obstacles. This assumption is at the heart of our approach, and our proofs no longer work without it. Notice however that such an assumption was not needed for the lower obstacles considered in \cite{mpz}. This is the first manifestation of an effect that we will highlight throughout the paper, namely that there is absolutely no symmetry between lower and upper obstacles in the second-order framework. 
\end{Remark}

\vspace{0.35em}
\begin{Remark}
The decomposition \reff{eq:decomp} is not restrictive. Indeed, with the integrability assumption (iv), we know that for each $\mathbb P\in\mathcal P^\kappa_H$, there exists a $\mathbb P$-martingale $M^\mathbb P$ and a bounded variation process $A^\mathbb P$ such that
$$S_t=S_0+M_t^\mathbb P+A_t^\mathbb P,\ \mathbb P-a.s.$$

Then, using the martingale representation theorem, there exists some $P_t^\mathbb P\in\mathbb H^2(\mathbb P)$ such that
$$M_t^\mathbb P=\int_0^tP_s^\mathbb PdB_s.$$

Then, since $S$ is c\`adl\`ag, by Karandikar \cite{kar}, we can aggregate the family $(P^\mathbb P)_{\mathbb P\in\mathcal P^\kappa_H}$ into a universal process $P$, which gives us the decomposition \reff{eq:decomp}.
\end{Remark}

Next, we also consider a lower obstacle $L$ which will be assumed to verify
\begin{Assumption}\label{assump.h3ref}
$\rm{(i)}$ $L$ is a $\mathbb F$-progressively measurable c\`adl\`ag process.

\vspace{0.35em}
$\rm{(ii)}$ $L$ is uniformly continuous in $\omega$ in the sense that for all $t$ and for some modulus of continuity $\rho$
$$\abs{L_t(\omega)-L_t(\widetilde{\omega})}\leq\rho\left(\No{\omega-\widetilde{\omega}}_t\right),\text{ }\forall \text{ }(\omega,\widetilde{\omega})\in\Omega^2.$$
\item[$\rm{(iii)}$] For all $t\in [0,T]$, we have
$$L_t<S_t\text{ and }L_{t^-}<S_{t^-},\ \mathcal P^\kappa_H-q.s.$$
\item[$\rm{(iv)}$] We have the following integrability condition
\begin{align}
\varphi^{2,\kappa}_H&:=\underset{\mathbb P\in\mathcal P^\kappa_H}{\sup}\mathbb E^{\mathbb P}\left[\underset{0\leq t\leq T}{\esup}^{\mathbb P}\left(\mathbb E_t^{H,\mathbb P}\left[\left(\underset{0\leq s\leq T}{\sup}\left(L_s\right)^{+}\right)^{\kappa}\right]\right)^{\frac{2}{\kappa}}\right]<+\infty.
\end{align}
\end{Assumption}

\begin{Remark}
Unlike for $S$, we did not assume here that $L$ was a semimartingale, and we cannot interchange the roles of $S$ and $L$, that is to say that the wellposedness results do not hold if we assume that $L$ is a semimartingale instead of $S$.
\end{Remark}

We shall consider the following second order doubly reflected BSDE ($2$DRBSDE for short) with upper obstacle $S$ and lower obstacle $L$
\begin{equation}
Y_t=\xi +\int_t^T\widehat{F}_s(Y_s,Z_s)ds -\int_t^TZ_sdB_s + V_T-V_t, \text{ } 0\leq t\leq T, \text{  } \mathcal P_H^\kappa-q.s.
\label{2dbsderef}
\end{equation}

In order to give the definition of the 2DRBSDE, we first need to introduce the corresponding standard doubly reflected BSDEs. Hence, for any $\mathbb{P}\in\mathcal{P}^\kappa_H$, $\mathbb{F}$-stopping time $\tau$, and $\mathcal{F}_\tau$-measurable random variable $\xi \in \mathbb{L}^2(\mathbb P)$, let $$(y^\mathbb{P},z^\mathbb{P},k^{\mathbb{P},+},k^{\mathbb P,-}):=(y^\mathbb{P}(\tau,\xi),z^\mathbb{P}(\tau,\xi),k^{\mathbb{P},+}(\tau,\xi),k^{\mathbb{P},-}(\tau,\xi)),$$ denote the unique solution to the following standard DRBSDE with upper obstacle $S$ and lower obstacle $L$ (existence and uniqueness have been proved under these assumptions in \cite{crep} among others)
\begin{equation}
\label{drbsde}
\begin{cases}
y_t^\mathbb P=\xi+\int_t^\tau\widehat F_s(y_s^\mathbb P,z_s^\mathbb P)ds-\int_t^\tau z_s^\mathbb PdB_s+k_{\tau}^{\mathbb{P},-}-k_{t}^{\mathbb{P},-}-k_{\tau}^{\mathbb{P},+}+k_{t}^{\mathbb{P},+}, \text{ }0\leq t\leq \tau,\text{ } \mathbb P-a.s.\\[0.8em]
L_t\leq y_t^\mathbb P\leq S_t, \text{ }\mathbb P-a.s.\\[0.8em]
\int_0^t\left(y^\mathbb P_{s^-}-L_{s^-}\right)dk^{\mathbb P,-}_s=\int_0^t\left(S_{s^-}-y^\mathbb P_{s^-}\right)dk^{\mathbb P,+}_s=0, \text{ }\mathbb P-a.s., \text{ }\forall t\in[0,T].
\end{cases}
\end{equation}

\begin{Remark}\label{neveract}
Notice that the assumption that $L_t<S_t$ and $L_{t^-}<S_{t^-}$ implies that the non-decreasing processes $k^{\mathbb P,+}$ and $k^{\mathbb P,-}$ never act at the same time. This will be important later. This hypothesis is already present in \cite{HH05} and \cite{HH06}.
\end{Remark}

\vspace{0.35em}
Everything is now ready for the
\begin{Definition}\label{death2}
We say $(Y,Z)\in \mathbb D^{2,\kappa}_H\times\mathbb H^{2,\kappa}_H$ is a solution to the 2DRBSDE \reff{2dbsderef} if

\begin{itemize}
\item[$\bullet$] $Y_T=\xi$, $\mathcal P_H^\kappa-q.s$.

\vspace{0.35em}
\item[$\bullet$] $\forall \mathbb P \in \mathcal P_H^\kappa$, the process $V^{\mathbb P}$ defined below has paths of bounded variation $\mathbb P-a.s.$
\begin{equation}
V_t^{\mathbb P}:=Y_0-Y_t - \int_0^t\widehat{F}_s(Y_s,Z_s)ds+\int_0^tZ_sdB_s, \text{ } 0\leq t\leq T, \text{  } \mathbb P-a.s.
\label{2dbsde.kref}
\end{equation}

\item[$\bullet$] We have the following minimum condition for $0\leq t\leq T$
\begin{equation}
V_t^\mathbb P+k_t^{\mathbb P,+}-k_t^{\mathbb P,-}=\underset{ \mathbb{P}^{'} \in \mathcal{P}_H(t^+,\mathbb{P}) }{ \einf^{\mathbb P} }\mathbb{E}_t^{\mathbb P^{'}}\left[V_T^{\mathbb P^{'}}+k_T^{\mathbb P^{'},+}-k_T^{\mathbb P^{'},-}\right], \text{  } \mathbb P-a.s., \text{ } \forall \mathbb P \in \mathcal P_H^\kappa.
\label{2dbsde.minref}
\end{equation}
\item[$\bullet$] $L_t\leq Y_t\leq S_t$, $\mathcal P_H^\kappa-q.s.$
\end{itemize}

Moreover, if there exists an aggregator for the family $(V^\mathbb P)_{\mathbb P\in\mathcal P^\kappa_H}$, that is to say a progressively measurable process $V$ such that for all $\mathbb P\in\mathcal P^\kappa_H$, $$V_t=V_t^\mathbb P,\ t\in[0,T],\ \mathbb P-a.s.,$$
then we say that $(Y,Z,V)$ is a solution to the 2DRBSDE \reff{2dbsderef}.

\end{Definition}

%

\vspace{0.35em}
\begin{Remark}\label{zhang}
The Definition \ref{death2} differs from the rest of the 2BSDE literature. Indeed, unlike in \cite{stz} for instance, the process $V^\mathbb P$ that we add in the definition of the 2BSDE is no longer non-decreasing, but is only assumed to have finite variation. This is mainly due to two competing effects. On the one hand, exactly as with standard RBSDEs with an upper obstacle, a non-increasing process has to be added to the solution in order to maintain it below the upper obstacle. But in the 2BSDE framework, another non-decreasing process also has to be added in order to push the process $Y$ to stay above all the $y^\mathbb P$ (as is shown by the representation formula proved below in Theorem \ref{uniqueref}). This emphasizes once more that in the second-order framework, which is fundamentally non-linear, there is no longer any symmetry between a reflected 2BSDE with an upper or a lower obstacle. Notice that this was to be expected, since 2BSDEs are a natural generalization of the $G$-expectation introduced by Peng \cite{peng}, which is an example of sublinear (and thus non-linear) expectation. We also would like to refer the reader to the recent paper by Pham and Zhang \cite{phamz}, whose problematics are strongly connected to ours. They study some norm estimates for semimartingales in the context of linear and sublinear expectations, and point out that there is a fundamental difference between non-linear submartingales and supermartingales (see their section $4.3$). Translated in our framework, and using the intuition from the classical RBSDE theory, when the generator is equal to $0$, a 2RBSDE with a lower obstacle should be a non-linear supermartingale, while a 2RBSDE with an upper obstacle should be a non-linear submartingale. In this sense, our results are a first step in the direction of the conjecture in section $4.3$ of \cite{phamz}.
\end{Remark}

\subsection{DRBSDEs as a special case of 2DRBSDEs}\label{connection}

In this subsection, we show how we can recover the usual theory. If $H$ is linear in $\gamma$, that is to say
$$H_t(y,z,\gamma):=\frac12\Tr{a_t^0\gamma}-f_t(y,z),$$
where $a^0:[0,T]\times\Omega\rightarrow \mathbb S_d^{>0}$ is $\mathbb F$-progressively measurable and has uniform upper and lower bounds, then as in \cite{stz}, we no longer need to assume any uniform continuity in $\omega$ in this case. Besides, the domain of $F$ is restricted to $a^0$ and we have
$$\widehat F_t(y,z)=f_t(y,z).$$

If we further assume that there exists some $\mathbb P\in \overline{\mathcal P}_S$ such that $\widehat a$ and $a^0$ coincide $\mathbb P-a.s.$ and $\mathbb E^\mathbb P\left[\int_0^T\abs{f_t(0,0)}^2dt\right]<+\infty$, then $\mathcal P^\kappa_H=\left\{\mathbb P\right\}$.

\vspace{0.35em}
Then, we know that $V^\mathbb P+k_t^{\mathbb P,+}-k_t^{\mathbb P,-}$ is a $\mathbb P$-martingale with finite variation. Since $\mathbb P$ satisfy the martingale representation property, this martingale is also continuous, and therefore it is null. Thus we have
$$0=k_t^{\mathbb P,+}-k_t^{\mathbb P,-}+V^\mathbb P,\text{ }\mathbb P-a.s.,$$
and the $2$DRBSDE is equivalent to a standard DRBSDE. In particular, we see that $V^\mathbb P$ now becomes a finite variation process which decreases only when $Y_{t^-}=S_{t^-}$ and increases only when $Y_{t^-}=L_{t^-}$. This implies
that $V^\mathbb P$ satisfies the usual Skorohod conditions. We would like to emphasize this fact here, since it will be useful later on to have a deeper understanding of the structure of the processes $\left\{V^\mathbb P\right\}_{\mathbb P\in\mathcal P^\kappa_H}$.

\section{Uniqueness, estimates and representations}\label{section.2}

\subsection{A representation inspired by stochastic control}
We have similarly as in Theorem $4.4$ of \cite{stz}

\begin{Theorem}\label{uniqueref}
Let Assumption \ref{assump.href} hold. Assume $\xi \in \mathbb{L}^{2,\kappa}_H$ and that $(Y,Z)$ is a solution to the $2$DRBSDE \reff{2dbsderef}. Then, for any $\mathbb{P}\in\mathcal{P}^\kappa_H$ and $0\leq t_1< t_2\leq T$,
\begin{align}
\label{representationref}
Y_{t_1}&=\underset{\mathbb{P}^{'}\in\mathcal{P}^\kappa_H(t_1^+,\mathbb{P})}{\esup^\mathbb{P}}y_{t_1}^{\mathbb{P}^{'}}(t_2,Y_{t_2}), \text{ }\mathbb{P}-a.s.
\end{align}

Consequently, the $2$DRBSDE \reff{2dbsderef} has at most one solution in $ \mathbb D^{2,\kappa}_H\times\mathbb H^{2,\kappa}_H$.
\end{Theorem}

\proof
The proof is exactly the same as the proof of Theorem $3.1$ in \cite{mpz}, so we will only sketch it. First, from the minimal condition \reff{2dbsde.minref}, we deduce that for any $\mathbb P\in\mathcal P^\kappa_H$, the process $V^\mathbb P+k_t^{\mathbb P,+}-k_t^{\mathbb P,-}$ is a $\mathbb P$-submartingale. By the Doob-Meyer decomposition and the martingale representation property, its martingale part is continuous with finite variation and therefore null. Hence $V^\mathbb P+k_t^{\mathbb P,+}-k_t^{\mathbb P,-}$ is a non-decreasing process. Then, the inequality
$$Y_{t_1}\geq\underset{\mathbb{P}^{'}\in\mathcal{P}^\kappa_H(t_1^+,\mathbb{P})}{\esup^\mathbb{P}}y_{t_1}^{\mathbb{P}^{'}}(t_2,Y_{t_2}), \text{ }\mathbb{P}-a.s.,$$
is a simple consequence of a classical comparison Theorem. The reverse inequality is then obtained by standard linearization techniques using the Lipschitz properties of $F$, see \cite{mpz} for the details.
\ep

\begin{Remark}\label{rem.min}
Let us now justify the minimum condition \reff{2dbsde.minref}. Assume for the sake of clarity that the generator $\widehat F$ is equal to $0$. By the above Theorem, we know that if there exists a solution to the $2$DRBSDE \reff{2dbsderef}, then the process $Y$ has to satisfy the representation \reff{representationref}. Therefore, we have a natural candidate for a possible solution of the $2$DRBSDE. Now, assume that we could construct such a process $Y$ satisfying the representation \reff{representationref} and which has the decomposition \reff{2dbsderef}. Then, taking conditional expectations in $Y-y^\mathbb P$, we end up with exactly the minimum condition \reff{2dbsde.minref}.
\end{Remark}

Finally, the following comparison Theorem follows easily from the classical one for DRBSDEs (see for instance \cite{lx2}) and the representation \reff{representationref}.

\begin{Theorem}
Let $(Y,Z)$ and $(Y',Z')$ {\rm (}resp. $(y^\mathbb P,z^\mathbb P,k^{+,\mathbb P},k^{-,\mathbb P})$ and $(y'^{\mathbb P},z'^{\mathbb P},k'^{+,\mathbb P},k'^{-,\mathbb P})${\rm )} be the solutions of the $2$DRBSDEs {\rm (}resp. DRBSDEs{\rm )} with terminal conditions $\xi$ and $\xi^{'}$, upper obstacles $S$ and $S^{'}$, lower obstacles $L$ and $L^{'}$ and generators $\widehat F$ and $\widehat F^{'}$ respectively. Assume that they both verify Assumptions \ref{assump.href}, that $\mathcal P^\kappa_H\subset\mathcal P^\kappa_{H^{'}}$ and that we have $\mathcal P^\kappa_H-q.s.$
$$\xi\leq\xi^{'}, \ \widehat F_t(y'^{\mathbb P}_t,z'^{\mathbb P}_t)\leq \widehat F^{'}_t(y'^{\mathbb P}_t,z'^{\mathbb P}_t),\ L_t\leq L_t^{'}\text{ and } S_t\geq S_t^{'}.$$
Then $Y\leq Y'$, $\mathcal P^\kappa_H-q.s.$
\end{Theorem}

\begin{Remark}
Unlike in the classical framework, even if the upper obstacles $S$ and $S^{'}$ and the lower obstacles $L$ and $L^{'}$ are identical, we cannot compare the processes $V^\mathbb P$ and $V'^{\mathbb P}$. This is due to the fact that these processes are not assumed to satisfy a Skorohod-type condition. This point was already mentioned in \cite{mpz}.
\end{Remark}

\subsection{A priori estimates}

We will now try to obtain a priori estimates for the 2DRBSDEs. We emphasize immediately that the fact that the process $V^\mathbb P$ are only of finite variation makes the task a lot more difficult than in \cite{mpz}. Indeed, we are now in a case which shares some similarities with standard doubly reflected BSDEs for which it is known that a priori estimates cannot be obtained without some regularity assumptions on the obstacles (for instance if one of them is a semimartingale). We assumed here that $S$ was a semimartingale, a property which will be at the heart of our proofs. Nonetheless, even before this, we need to understand the fine structure of the processes $V^\mathbb P$. This is the object of the following proposition.

\vspace{0.5em}
\begin{Proposition}\label{prop.imp}
Let Assumption \ref{assump.href} hold. Assume $\xi\in\mathbb L^{2,\kappa}_H$ and $(Y,Z)\in \mathbb D^{2,\kappa}_H\times\mathbb H^{2,\kappa}_H$ is a solution to the 2DRBSDE \reff{2dbsderef}. Let $\left\{(y^\mathbb P,z^\mathbb P,k^{+,\mathbb P},k^{-,\P})\right\}_{\mathbb P\in\mathcal P^\kappa_H}$ be the solutions of the corresponding DRBSDEs \reff{drbsde}. Then we have the following results for all $t\in [0,T]$ and for all $\mathbb P\in\mathcal P^\kappa_H$

\vspace{0.55em}
$\rm{(i)}$ $V_t^{\mathbb P,+}:=\int_0^t{\bf 1}_{y_{s^-}^\mathbb P<S_{s^-}}dV_s^\mathbb P$ is a non-decreasing process, $\mathbb P-a.s.$

\vspace{0.55em}
$\rm{(ii)}$ $V_t^{\mathbb P,-}:=\int_0^t{\bf 1}_{y_{s^-}^\mathbb P=S_{s^-}}dV_s^\mathbb P=-k_t^{\mathbb P,+},$ $\mathbb P-a.s.,$ and is therefore a non-increasing process.
\end{Proposition}

\vspace{0.5em}
\proof

\vspace{0.5em}
Let us fix a given $\mathbb P\in\mathcal P^\kappa_H$.

\vspace{0.35em}
\rm{(i)} Let $\tau_1$ and $\tau_2$ be two $\mathbb F$-stopping times such that for all $t\in[\tau_1,\tau_2)$, $y_{t^-}^\mathbb P<S_{t^-}$, $\mathbb P-a.s.$

\vspace{0.55em}
Then, we know from the usual Skorohod condition that the process $k^{\mathbb P,+}$ does not increase between $\tau_1$ and $\tau_2$. Now, we remind the reader that we showed in the proof of Theorem \ref{uniqueref}, that the process $V^\mathbb P+k^{+,\mathbb P}-k^{-,\mathbb P}$ is always non-decreasing. This necessarily implies that $V^\mathbb P$ must be non-decreasing between $\tau_1$ and $\tau_2$. Hence the first result.

\vspace{0.55em}
\rm{(ii)} Let now $\tau_1$ and $\tau_2$ be two $\mathbb F$-stopping times such that for all $t\in[\tau_1,\tau_2)$, $y_{t^-}^\mathbb P=S_{t^-}$, $\mathbb P-a.s.$

\vspace{0.55em}
First, since the two obstacles are separated, we necessarily have $y_{t^-}^\P>L_{t^-}$, $\P-a.s.$ for every $t\in[\tau_1,\tau_2)$, which in turn implies that $k^{\P,-}$ does not increase. Next, by the representation formula \reff{representationref}, we necessarily have $Y_{t^-}\geq y^\mathbb P_{t^-}$, $\mathbb P-a.s.$ for all $t$. Moreover, since we also have $Y_t\leq S_t$ by definition, this implies, since all the processes here are c\`adl\`ag, that we must have
$$Y_{t^-}=y^\mathbb P_{t^-}=S_{t^-},\ t\in[\tau_1,\tau_2),\ \mathbb P-a.s.$$

Using the fact that $Y$ and $y^\mathbb P$ solve respectively a 2BSDE and a BSDE, we also have $\mathbb P-a.s.$
\begin{align*}
&S_{t^{-}}+\Delta Y_t=Y_t=Y_{u}+\int_t^{u}\widehat F_s(Y_s,Z_s)ds-\int_t^{u}Z_sdB_s+V_{u}^\mathbb P-V_t^{\mathbb P},\ \tau_1\leq t\leq u< \tau_2\\
&S_{t^{-}}+\Delta y^\mathbb P_t=y_t^\mathbb P=y^\mathbb P_{u}+\int_t^{u}\widehat F_s(y^\mathbb P_s,z^\mathbb P_s)ds-\int_t^{u}z^\mathbb P_sdB_s-k_{u}^{\mathbb P,+}+k_t^{\mathbb P,+},\ \tau_1\leq t\leq u< \tau_2.
\end{align*}

Identifying the martingale parts above, we obtain that $Z_s=z^\mathbb P_s$, $ds\times\mathbb P-a.e.$ Then, identifying the finite variation parts, we have
\begin{align}\label{pouic}
&Y_{u}-\Delta Y_t+\int_t^u\widehat F_s(Y_s,Z_s)ds+V_{u}^\mathbb P-V_t^{\mathbb P}=y^\mathbb P_{u}-\Delta y^\mathbb P_t+\int_t^u\widehat F_s(y^\mathbb P_s,z^\mathbb P_s)ds-k_{u}^{\mathbb P,+}+k_t^{\mathbb P,+}.
\end{align}

Now, we clearly have
$$\int_t^u\widehat F_s(Y_s,Z_s)ds=\int_t^u\widehat F_s(y^\mathbb P_s,z^\mathbb P_s)ds,$$
since $Z_s=z^\mathbb P_s$, $dt\times\mathbb P-a.e.$  and $Y_{s^-}=y^\mathbb P_{s^-}=S_{s^-}$ for all $s\in [t,u]$. Moreover, since $Y_{s^-}=y^\mathbb P_{s^-}=S_{s^-}$ for all $s\in [t,u]$ and since all the processes are c\`adl\`ag, the jumps of $Y$ and $y^\mathbb P$ are equal to the jumps of $S$. Therefore, \reff{pouic} can be rewritten
$$V_{u}^\mathbb P-V_t^{\mathbb P}=-k_{u}^{\mathbb P,+}+k_t^{\mathbb P,+},$$
which is the desired result.
\ep

\vspace{0.35em}
The above Proposition is crucial for us. Indeed, we have actually shown that
$$V_t^\mathbb P=V_t^{\mathbb P,+}-k_t^{\mathbb P,+},\ \mathbb P-a.s.,$$
where $V^{\mathbb P,+}$ and $k^{\mathbb P,+}$ are two non-decreasing processes which never act at the same time. Hence, we have obtained a Jordan decomposition for $V^\P$. Moreover, we can easily obtain a priori estimates for $k^{\mathbb P,+}$ by using the fact that it is part of the solution of the DRBSDE \reff{drbsde}. Notice that these estimates hold only because we assumed that the corresponding upper obstacle $S$ was a semimartingale. This is this decomposition which will allow us to obtain the estimates.

\begin{Remark}
The above result is enough for us to obtain the desired a priori estimates. However, we can go further into the structure of the bounded variation processes $V^\mathbb P$. Indeed, arguing as in Proposition $3.2$ in \cite{mpz}, we could also show that
$${\bf 1}_{Y_{t^-}=L_{t^-}}dV_t^{\mathbb P}= {\bf 1}_{Y_{t^-}=L_{t^-}}dk_t^{\mathbb P,-}.$$

Notice however that we, a priori, cannot say anything about $V^\mathbb P$ when $L_{t^-}=y_{t^-}^{\mathbb P}<Y_{t^-}$, even though we showed that it could be known explicitly when $S_{t^-}=y_{t^-}^{\mathbb P}$. This emphasizes once more the fact that the upper and the lower obstacle in our context do not play a symmetric role.
\end{Remark}

\vspace{0.35em}
We can now prove the following theorem.
\begin{Theorem}\label{destimatesref}
Let Assumptions \ref{assump.href}, \ref{assump.obs} and \ref{assump.h3ref} hold. Assume $\xi\in\mathbb L^{2,\kappa}_H$ and $(Y,Z)\in \mathbb D^{2,\kappa}_H\times\mathbb H^{2,\kappa}_H$ is a solution to the 2DRBSDE \reff{2dbsderef}. Let $\left\{(y^\mathbb P,z^\mathbb P,k^{\mathbb P,+},k^{\mathbb P,-})\right\}_{\mathbb P\in\mathcal P^\kappa_H}$ be the solutions of the corresponding DRBSDEs \reff{drbsde}. Then, there exists a constant $C_\kappa$ depending only on $\kappa$, $T$ and the Lipschitz constant of $\widehat F$ such that
\begin{align*}
&\No{Y}^2_{\mathbb D^{2,\kappa}_H}+\No{Z}^2_{\mathbb H^{2,\kappa}_H}+\underset{\mathbb P\in \mathcal P^\kappa_H}{\sup}\left\{\No{y^\mathbb P}^2_{\mathbb D^{2}(\mathbb P)}+\No{z^\mathbb P}^2_{\mathbb H^{2}(\mathbb P)}+\displaystyle\mathbb E^\mathbb P\left[{\rm{Var}}_{0,T}\left(V^\mathbb P\right)^2+\left(k_T^{\mathbb P,+}\right)^2+\left(k_T^{\mathbb P,-}\right)^2\right]\right\}\\
&\leq C_\kappa\left(\No{\xi}^2_{\mathbb L^{2,\kappa}_H}+\phi^{2,\kappa}_H+\psi^{2,\kappa}_H+\varphi^{2,\kappa}_H+\zeta^{2,\kappa}_H\right)
\end{align*}
\end{Theorem}

\proof
First of all, since we assumed that $S$ was a semimartingale, we can argue as in \cite{crep} to obtain that
\begin{align*}
dk_t^{\mathbb P,+}&\leq \widehat F_t^+(S_t,P_t)dt+dA_t^{\mathbb P,+}\leq C\left(\abs{\widehat F_t^0}+\abs{S_t}+\abs{\widehat a_t^{1/2}P_t}\right)dt+dA_t^{\mathbb P,+}.
\end{align*}

Hence,
\begin{align*}
\mathbb E_t^\mathbb P\left[\left(k_T^{\mathbb P,+}\right)^\kappa\right] &\leq C_\kappa\mathbb E_t^\mathbb P\left[ \int_t^T\abs{\widehat F_s^0}^\kappa ds+\left(\int_t^T|\widehat a_s^{\frac12}P_s|^2ds\right)^{\kappa/2}+\underset{t\leq s\leq T}{\sup}\abs{S_s}^\kappa+\left(A_T^{\mathbb P,-}\right)^\kappa\right]\\
&\leq C_\kappa\left(\left(\zeta^{2,\kappa}_H\right)^{1/2}+\mathbb E_t^\mathbb P\left[\int_t^T\abs{\widehat F_s^0}^\kappa ds+\underset{t\leq s\leq T}{\sup}\abs{S_s}^\kappa\right]\right).
\end{align*}

Let us now define
$$\tilde \xi:=\xi-k_T^{\mathbb P,+},\ \tilde y^\mathbb P=y^\mathbb P-k^{\mathbb P,+},\ \tilde F_t(y,z):=\widehat F_t\left(y+k^{\mathbb P,+}_t,z\right).$$

Then, it is easy to see that $(\tilde y^\mathbb P,z^\mathbb P, k^{\mathbb P,+})$ is the solution of the lower reflected BSDE with terminal condition $\tilde \xi$, generator $\tilde F$ and obstacle $L-k^{\mathbb P,+}$. We can then once again apply Lemma $2$ in \cite{ham} to obtain that there exists a constant $C_\kappa$ depending only on $\kappa$, $T$ and the Lipschitz constant of $\widehat F$, such that for all $\mathbb P$
\begin{align}
\label{estimrefd}
\nonumber\abs{y_t^\mathbb P}&\leq C_\kappa\mathbb E_t^\mathbb P\left[\abs{\tilde\xi}^\kappa+\int_t^T\abs{\tilde F^0_s}^\kappa ds+\underset{t\leq s\leq T}{\sup}\left(\left(L_s-k_s^{\mathbb P,+}\right)^+\right)^\kappa\right]\\
\nonumber&\leq C_\kappa\mathbb E_t^\mathbb P\left[\abs{\xi}^\kappa+\int_t^T\abs{\widehat F^0_s}^\kappa ds+\underset{t\leq s\leq T}{\sup}\left(L_s^+\right)^\kappa+\left(k_T^{\mathbb P,+}\right)^\kappa\right]\\
&\leq C_\kappa\left(\left(\zeta^{2,\kappa}_H\right)^{1/2}+\mathbb E_t^\mathbb P\left[\abs{\xi}^\kappa+\int_t^T\abs{\widehat F^0_s}^\kappa ds+\underset{t\leq s\leq T}{\sup}\abs{S_s}^\kappa+\underset{t\leq s\leq T}{\sup}\left(L_s^+\right)^\kappa\right]\right).
\end{align}

This immediately provides the estimate for $y^\mathbb P$. Now by definition of the norms, we obtain from \reff{estimrefd} and the representation formula \reff{representationref} that
\begin{equation}
\label{eq.yd}
\No{Y}_{\mathbb D^{2,\kappa}_H}^2\leq C_\kappa\left(\No{\xi}^2_{\mathbb L^{2,\kappa}_H}+\phi^{2,\kappa}_H+\psi^{2,\kappa}_H+\varphi^{2,\kappa}_H+\zeta^{2,\kappa}_H\right).
\end{equation}

Now apply It\^o's formula to $\abs{y^\mathbb P}^2$ under $\mathbb P$. We get as usual for every $\epsilon>0$
\begin{align}
\nonumber\mathbb E^\mathbb P\left[\int_0^T\abs{\widehat a_t^{1/2}z_t^\mathbb P}^2dt\right]&\leq C\mathbb E^\mathbb P\left[\abs{\xi}^2+\int_0^T\abs{y^\mathbb P_t}\left(\abs{\widehat F^0_t}+\abs{y_t^\mathbb P}+\abs{\widehat a_t^{1/2}z^\mathbb P_t}\right)dt\right]\\
&\nonumber\hspace{0.9em}+\mathbb E^\mathbb P\left[\int_0^T\abs{y_t^\mathbb P}d\left(k_t^{\mathbb P,+}+k_t^{\mathbb P,-}\right)\right]\\
\nonumber&\leq C\left(\No{\xi}_{\mathbb L^{2,\kappa}_H}+\mathbb E^\mathbb P\left[\underset{0\leq t\leq T}{\sup}\abs{y^\mathbb P_t}^2+\left(k_T^{\mathbb P,-}\right)^2+\left(\int_0^T\abs{\widehat F^0_t}dt\right)^2\right]\right)\\
&\hspace{0.9em}+\epsilon\mathbb E^\mathbb P\left[\int_0^T\abs{\widehat a_t^{1/2}z^\mathbb P_t}^2dt+\abs{k_T^{\mathbb P,+}}^2\right]+\frac{C^2}{\eps}\mathbb E^\mathbb P\left[\underset{0\leq t\leq T}{\sup}\abs{y^\mathbb P_t}^2\right].
\label{grahou2d}
\end{align}

Then by definition of the DRBSDE \reff{drbsde}, we easily have
\begin{equation}
\mathbb E^\mathbb P\left[\abs{k_T^{\mathbb P,+}}^2\right]\leq C_0\mathbb E^\mathbb P\left[\abs{\xi}^2+\underset{0\leq t\leq T}{\sup}\abs{y^\mathbb P_t}^2+\left(k_T^{\mathbb P,-}\right)^2+\int_0^T\abs{\widehat a_t^{\frac12}z^\mathbb P_t}^2dt+\left(\int_0^T\abs{\widehat F^0_t}dt\right)^2\right],
\label{eq.kk2d}
\end{equation}
for some constant $C_0$, independent of $\epsilon$. Now set $\epsilon:=(2(1+C_0))^{-1}$ and plug \reff{eq.kk2d} in \reff{grahou2d}. We obtain from the estimates for $y^\mathbb P$ and $k^{\mathbb P,-}$
$$\underset{\mathbb P\in \mathcal P^\kappa_H}{\sup}\No{z^\mathbb P}_{\mathbb H^{2}(\mathbb P)}\leq C\left(\No{\xi}^2_{\mathbb L^{2,\kappa}_H}+\phi^{2,\kappa}_H+\psi^{2,\kappa}_H+\varphi^{2,\kappa}_H+\zeta^{2,\kappa}_H\right).$$

Then the estimate for $k^{\mathbb P,+}$ comes from \reff{eq.kk2d}. Now that we have obtained the desired estimates for $y^\mathbb P$ ,$z^\mathbb P$, $k^{\mathbb P,+}$, $k^{\mathbb P,-}$ and $Y$, we can proceed further. 

\vspace{0.35em}
Exactly as above, we apply It\^o's formula to $\abs{Y}^2$ under each $\mathbb P\in\mathcal P^\kappa_H$. We have once more for every $\epsilon>0$ and using Proposition \ref{prop.imp}
\begin{align}
\nonumber\mathbb E^\mathbb P\left[\int_0^T|\widehat a_t^{\frac12}Z_t|^2dt\right]&\leq C\mathbb E^\mathbb P\left[\abs{\xi}^2+\int_0^T\abs{Y_t}\left(\abs{\widehat F^0_t}+\abs{Y_t}+|\widehat a_t^{\frac12}Z_t|\right)dt\right]\\
&\nonumber\hspace{0.9em}+\mathbb E^\mathbb P\left[\int_0^TY_tdV_t^{\mathbb P,+}-\int_0^TY_tdk_t^{\mathbb P,+}\right]\\
\nonumber&\leq C\left(\No{\xi}_{\mathbb L^{2,\kappa}_H}+\mathbb E^\mathbb P\left[\underset{0\leq t\leq T}{\sup}\abs{Y_t}^2+\left(\int_0^T\abs{\widehat F^0_t}dt\right)^2\right]\right)\\
&\hspace{0.9em}+\epsilon\mathbb E^\mathbb P\left[\int_0^T|\widehat a_t^{\frac12}Z_t|^2dt+\abs{k_T^{\mathbb P,+}}^2+\abs{V_T^{\mathbb P,+}}^2+\frac{C^2}{\eps}\underset{0\leq t\leq T}{\sup}\abs{Y_t}^2\right].
\label{grahoud}
\end{align}

\vspace{0.35em}
Then by definition of our $2$DRBSDE, we easily have
\begin{equation}
\mathbb E^\mathbb P\left[\abs{V_T^{\mathbb P,+}}^2\right]\leq C_0\mathbb E^\mathbb P\left[\abs{\xi}^2+\underset{\scriptstyle0\leq t\leq T}{\sup}\abs{Y_t}^2+\int_0^T|\widehat a_t^{\frac12}Z_t|^2dt+\abs{k_T^{\mathbb P,+}}^2+\left(\int_0^T|\widehat F^0_t|dt\right)^2\right],
\label{eq.kkd}
\end{equation}
for some constant $C_0$, independent of $\epsilon$.

\vspace{0.35em}
Now set $\epsilon:=(2(1+C_0))^{-1}$ and plug \reff{eq.kkd} in \reff{grahoud}. One then gets
$$\mathbb E^\mathbb P\left[\int_0^T\abs{\widehat a_t^{1/2}Z_t}^2dt\right]\leq C\mathbb E^\mathbb P\left[\abs{\xi}^2+\underset{0\leq t\leq T}{\sup}\abs{Y_t}^2+\abs{k_T^{\mathbb P,+}}^2+\left(\int_0^T\abs{\widehat F^0_t}dt\right)^2\right].$$

From this and the estimates for $Y$ and $k^{\mathbb P,+}$, we immediately obtain
$$\No{Z}_{\mathbb H^{2,\kappa}_H}\leq C\left(\No{\xi}^2_{\mathbb L^{2,\kappa}_H}+\phi^{2,\kappa}_H+\psi^{2,\kappa}_H+\varphi^{2,\kappa}_H+\zeta^{2,\kappa}_H\right).$$

Moreover, we deduce from \reff{eq.kkd} that
\begin{equation}
\underset{\mathbb P\in\mathcal P^\kappa_H}{\sup}\mathbb E^\mathbb P\left[\left(V^{\mathbb P,+}_T\right)^2\right]\leq C\left(\No{\xi}^2_{\mathbb L^{2,\kappa}_H}+\phi^{2,\kappa}_H+\psi^{2,\kappa}_H+\varphi^{2,\kappa}_H+\zeta^{2,\kappa}_H\right).
\label{eqqd}
\end{equation}

Finally, we have by definition of the total variation and the fact that the processes $V^{\mathbb P,+}$ and $k^{\mathbb P,+}$ are non-decreasing
\begin{align*}
\mathbb E^\mathbb P\left[\underset{0,T}{\rm{Var}}\left(V^\mathbb P\right)^2\right]&\leq C\mathbb E^\mathbb P\left[\underset{0,T}{\rm{Var}}\left(V^{\mathbb P,+}\right)^2+\underset{0,T}{\rm{Var}}\left(k^{\mathbb P,+}\right)^2\right]\\
&=C\mathbb E^\mathbb P\left[\left(V_T^{\mathbb P,+}\right)^2+\left(k_T^{\mathbb P,+}\right)^2\right]\\
&\leq C\left(\No{\xi}^2_{\mathbb L^{2,\kappa}_H}+\phi^{2,\kappa}_H+\psi^{2,\kappa}_H+\varphi^{2,\kappa}_H+\zeta^{2,\kappa}_H\right),
\end{align*}
where we used the estimate for $k^{\mathbb P,+}$ and \reff{eqqd} for the last inequality.
\ep

\vspace{0.35em}
\begin{Theorem}\label{estimates2d}
Let Assumptions \ref{assump.href}, \ref{assump.obs} and \ref{assump.h3ref} hold. For $i=1,2$, let $(Y^i,Z^i)$ be the solutions to the 2DRBSDE \reff{2dbsderef} with terminal condition $\xi^i$, upper obstacle $S$ and lower obstacle $L$. Then, there exists a constant $C_\kappa$ depending only on $\kappa$, $T$ and the Lipschitz constant of $F$ such that
\begin{align*}
&\No{Y^1-Y^2}_{\mathbb D^{2,\kappa}_H}\leq C\No{\xi^1-\xi^2}_{\mathbb L^{2,\kappa}_H}\\
&\No{Z^1-Z^2}^2_{\mathbb H^{2,\kappa}_H}+\underset{\mathbb P\in \mathcal P^\kappa_H}{\sup}\mathbb E^\mathbb P\left[\underset{0\leq t\leq T}{\sup}\abs{V_t^{\mathbb P,+,1}-V_t^{\mathbb P,+,2}}^2+\underset{0\leq t\leq T}{\sup}\abs{V_t^{\mathbb P,-,1}-V_t^{\mathbb P,-,2}}^2\right]\\
&\leq C\No{\xi^1-\xi^2}_{\mathbb L^{2,\kappa}_H}\left(\No{\xi^1}_{\mathbb L^{2,\kappa}_H}+\No{\xi^1}_{\mathbb L^{2,\kappa}_H}+(\phi^{2,\kappa}_H)^{1/2}+(\psi^{2,\kappa}_H)^{1/2}+(\varphi^{2,\kappa}_H)^{1/2}+(\zeta^{2,\kappa}_H)^{1/2}\right).
\end{align*}
\end{Theorem}

\begin{Remark}
We emphasize that in Theorem \ref{estimates2d}, we control the norm of both $V_t^{\mathbb P,+,1}-V_t^{\mathbb P,+,2}$ and $V_t^{\mathbb P,-,1}-V_t^{\mathbb P,-,2}$. This is crucial in our main existence Theorem \ref{mainrefd}.
\end{Remark}

\proof
As in the previous Proposition, we can follow the proof of Lemma $3$ in \cite{ham}, to obtain that there exists a constant $C_\kappa$ depending only on $\kappa$, $T$ and the Lipschitz constant of $F$, such that for all $\mathbb P$
\begin{equation}
\label{estim2d}
\abs{y_t^{\mathbb P,1}-y_t^{\mathbb P,2}}\leq C_\kappa\mathbb E_t^\mathbb P\left[\abs{\xi^1-\xi^2}^\kappa\right].
\end{equation}

Now by definition of the norms, we get from \reff{estim2d} and the representation formula \reff{representationref} that
\begin{equation}
\label{eq.y2d}
\No{Y^1-Y^2}_{\mathbb D^{2,\kappa}_H}^2\leq C_\kappa\No{\xi^1-\xi^2}^2_{\mathbb L^{2,\kappa}_H}.
\end{equation}

Next, the estimate for $V_t^{\mathbb P,-,1}-V_t^{\mathbb P,-,2}$ is immediate from the usual estimates for DRBSDEs (see for instance Theorem $3.2$ in \cite{crep}), since we actually have from Proposition \ref{prop.imp}
$$V_t^{\mathbb P,-,1}-V_t^{\mathbb P,-,2}=k^{\P,+,2}-k_t^{\P,+,1}.$$

Applying It\^o's formula to $\abs{Y^1-Y^2}^2$, under each $\mathbb P\in\mathcal P^\kappa_H$, leads to
\begin{align*}
\nonumber\mathbb E^\mathbb P\left[\int_0^T\abs{\widehat a_t^{\frac12}(Z_t^1-Z_t^2)}^2dt\right]&\leq C\mathbb E^\mathbb P\left[\abs{\xi^1-\xi^2}^2\right]+\mathbb E^\mathbb P\left[\int_0^T(Y_t^1-Y_t^2)d(V_t^{\mathbb P,1}-V_t^{\mathbb P,2})\right]\\
\nonumber&\hspace{0.9em}+C\mathbb E^\mathbb P\left[\int_0^T\abs{Y_t^1-Y_t^2}\left(\abs{Y_t^1-Y_t^2}+  |\widehat a_t^{\frac12}(Z_t^1-Z_t^2)|\right)dt\right]\\[0.8em]
\nonumber&\leq C\left(\No{\xi^1-\xi^2}_{\mathbb L^{2,\kappa}_H}^2+\No{Y^1-Y^2}^2_{\mathbb D^{2,\kappa}_H}\right)\\
\nonumber&\hspace{0.9em}+\frac12\mathbb E^\mathbb P\left[\int_0^T\abs{\widehat a_t^{1/2}(Z_t^1-Z_t^2)}^2dt\right]\\
&\hspace{0.9em}+C\No{Y^1-Y^2}_{\mathbb D^{2,\kappa}_H}\left(\mathbb E^\mathbb P\left[\sum_{i=1}^2\rm{Var}_{0,T}\left(V^{\mathbb P,i}\right)^2\right]\right)^{\frac12}.
\end{align*}

The estimate for $(Z^1-Z^2)$ is now obvious from the above inequality and the estimates of Proposition \ref{destimatesref}. Finally, we have by definition for any $t\in[0,T]$
\begin{align*}
\mathbb E^\mathbb P\left[\underset{0\leq t\leq T}{\sup}\abs{V_t^{\mathbb P,+,1}-V_t^{\mathbb P,+,2}}^2\right]&\leq \mathbb E^\mathbb P\left[\abs{\xi^1-\xi^2}^2+\underset{0\leq t\leq T}{\sup}\left\{\abs{Y_t^1-Y_t^2}^2+\abs{\int_t^T(Z_s^1-Z_s^2)dB_s}^2\right\}\right]\\
&\hspace{0.9em}+\mathbb E^\mathbb P\left[\int_0^T\abs{\widehat F_s(Y_s^1,Z_s^1)-\widehat F_s(Y_s^2,Z_s^2)}^2ds+\underset{0\leq t\leq T}{\sup}\abs{V_t^{\mathbb P,-,1}-V_t^{\mathbb P,-,2}}^2\right]\\
&\leq C\left(\No{\xi^1-\xi^2}_{\mathbb L^{2,\kappa}_H}^2+\No{Y^1-Y^2}_{\mathbb D^{2,\kappa}_H}^2+\No{Z^1-Z^2}_{\mathbb H^{2,\kappa}_H}^2\right),
\end{align*}
where we used the BDG inequality for the last step.

\vspace{0.35em}
By all the previous estimates, this finishes the proof.
\ep

\subsection{Some properties of the solution}

Now that we have proved the representation \reff{representationref} and the a priori estimates of Theorems \ref{destimatesref} and \ref{estimates2d}, we can show, as in the classical framework, that the solution $Y$ of the $2$DRBSDE is linked to some kind of Dynkin game. We emphasize that such a connection with games was already conjectured in \cite{mpz}. After that, Bayraktar and Yao \cite{bay} showed in a purely Markovian context, that the value function of stochastic zero-sum differential game could be linked to the notion of 2DRBSDEs, even though these objects were not precisely defined in the paper (see their section $5.2$). For any $t\in[0,T],$ denote $\mathcal T_{t,T}$ the set of $\F$-stopping times taking values in $[t,T]$.

\begin{Proposition}\label{dynkin}
Let $(Y,Z)$ be the solution to the above $2$DRBSDE \reff{2dbsderef}. For any $(\tau,\sigma)\in\mathcal T_{0,T}$, define
$$R_\tau^\sigma:=S_\tau1_{\tau<\sigma}+L_\sigma1_{\sigma\leq \tau, \sigma<T}+\xi1_{\tau\wedge\sigma=T}.$$
Then for each $t\in [0,T]$, for all $\mathbb P\in\mathcal P^\kappa_H$, we have $\mathbb P-a.s.$
\begin{align*}
Y_t&=\underset{\mathbb P^{'}\in\mathcal P^\kappa_H(t^+,\mathbb P)}{\esup^\mathbb P}\underset{\tau\in\mathcal T_{t,T}}{\einf}\ \underset{\sigma\in\mathcal T_{t,T}}{\esup}\text{ }\mathbb E_t^{\mathbb P^{'}}\left[\int_t^{\tau\wedge\sigma}\widehat F_s(y_s^{\mathbb P^{'}},z_s^{\mathbb P^{'}})ds +R_\tau^\sigma\right]\\
&=\underset{\mathbb P^{'}\in\mathcal P^\kappa_H(t^+,\mathbb P)}{\esup^\mathbb P}\underset{\sigma\in\mathcal T_{t,T}}{\esup}\ \underset{\tau\in\mathcal T_{t,T}}{\einf}\text{ }\mathbb E_t^{\mathbb P^{'}}\left[\int_t^{\tau\wedge\sigma}\widehat F_s(y_s^{\mathbb P^{'}},z_s^{\mathbb P^{'}})ds +R_\tau^\sigma\right].
\end{align*}

Moreover, for any $\gamma\in[0,1]$, we have $\mathbb P-a.s.$
\begin{align*}
Y_t&=\underset{\tau\in\mathcal T_{t,T}}{\einf}\ \underset{\sigma\in\mathcal T_{t,T}}{\esup}\text{ }\mathbb E_t^\mathbb P\left[\int_t^{\tau\wedge\sigma}\widehat F_s(
Y_s,Z_s)ds +K_{\tau\wedge\sigma}^{\mathbb P,\gamma}-K_t^{\mathbb P,\gamma}+R_\tau^\sigma\right]\\
&=\underset{\sigma\in\mathcal T_{t,T}}{\esup}\ \underset{\tau\in\mathcal T_{t,T}}{\einf}\text{ }\mathbb E_t^\mathbb P\left[\int_t^{\tau\wedge\sigma}\widehat F_s(
Y_s,Z_s)ds +K_{\tau\wedge\sigma}^{\mathbb P,\gamma}-K_t^{\mathbb P,\gamma}+R_\tau^\sigma\right],
\end{align*}
where $$K_t^{\mathbb P,\gamma}:=\gamma\int_0^t{\bf 1}_{y_{s^-}^\mathbb P<S_{s^-}}dV_s^\mathbb P+(1-\gamma)\int_0^t{\bf 1}_{Y_{s^-}>L_{s^-}}dV_s^\mathbb P.$$

Furthermore, for any $\mathbb P\in\mathcal P^\kappa_H$, the following stopping times are $\eps$-optimal
$$\tau_t^{\eps,\mathbb P}:=\inf\left\{s\geq t,\ y_s^\mathbb P\geq S_s-\eps,\ \mathbb P-a.s.\right\}\text{ and }\sigma_t^{\eps}:=\inf\left\{s\geq t,\ Y_s\leq L_s+\eps,\ \mathcal P^\kappa_H-q.s.\right\}.$$
\end{Proposition}

\begin{Remark}
Notice that the optimal stopping rules above are different in nature. Indeed, $\tau_t^{\eps,\mathbb P}$ depends explicitly on the probability measures $\mathbb P$, because it depends on the process $y^\mathbb P$, while $\sigma_t^\eps$ only depends on $Y$. This situation shed once more light on the complete absence of symmetry between lower and upper obstacles in the second-order framework
\end{Remark}

\begin{Remark}
The second result in the Proposition above may seem peculiar at first sight because of the degree of freedom introduced by the parameter $\gamma$. However, as shown in the proof below, we can find stopping times which are $\eps$-optimizer for the corresponding stochastic game, and which roughly correspond (as expected) to the first hitting times of the obstacles. Since the latter are completely separated, we know from Proposition \ref{prop.imp} that before hitting $S$, $$dV_t^\mathbb P={\bf 1}_{y_{t^-}^\mathbb P<S_{t^-}},$$ and that before hitting $L$, $$dV_t^\mathbb P={\bf 1}_{Y_{t^-}>L_{t^-}}.$$

Thanks to this result, it is easy to see that we can change the value of $\gamma$ as we see fit. In particular, if there is no upper obstacle, that is to say if $S=+\infty$, then taking $\gamma=0$, we recover the result of Proposition $3.1$ in \cite{mpz}.
\end{Remark}

\vspace{0.35em}
\proof
By Proposition $3.1$ in \cite{lx2}, we know that for all $\mathbb P\in\mathcal P^\kappa_H$, $\mathbb P-a.s.$
\begin{align*}
y_t^\mathbb P&=\underset{\tau\in\mathcal T_{t,T}}{\einf}\ \underset{\sigma\in\mathcal T_{t,T}}{\esup}\text{ }\mathbb E_t^{\mathbb P}\left[\int_t^{\tau\wedge\sigma}\widehat F_s(y_s^{\mathbb P},z_s^{\mathbb P})ds +R_\tau^\sigma\right]=\underset{\sigma\in\mathcal T_{t,T}}{\esup}\ \underset{\tau\in\mathcal T_{t,T}}{\einf}\text{ }\mathbb E_t^{\mathbb P}\left[\int_t^{\tau\wedge\sigma}\widehat F_s(y_s^{\mathbb P},z_s^{\mathbb P})ds +R_\tau^\sigma\right].
\end{align*}

Then the first equality is a simple consequence of the representation formula \reff{representationref}. For the second one, we proceed exactly as in the proof of Proposition $3.1$ in \cite{lx2}. Fix some $\mathbb P\in\mathcal P^\kappa_H$ and some $t\in[0,T]$ and some $\eps>0$. It is then easy to show that for any $s\in[t,\tau_t^{\eps,\mathbb P}]$, we have $y^\mathbb P_{s^-}<S_{s^-}$. In particular this implies that
$$dV_s^\mathbb P={\bf 1}_{y^\mathbb P_{s^-}<S_{s^-}}dV_s^\mathbb P,\ s\in[t,\tau_t^{\eps,\mathbb P}].$$

Let now $\sigma\in\mathcal T_{t,T}$. On the set $\{\tau^{\eps,\mathbb P}_t<\sigma\}$, we have
\begin{align*}
&\int_t^{\sigma\wedge\tau^{\eps,\mathbb P}_t}\widehat F_s(Y_s,Z_s)ds+R_{\tau^{\eps,\mathbb P}_t}^\sigma=\int_t^{\tau^{\eps,\mathbb P}_t}\widehat F_s(Y_s,Z_s)ds+S_{\tau^{\eps,\mathbb P}_t}\leq \int_t^{\tau^{\eps,\mathbb P}_t}\widehat F_s(Y_s,Z_s)ds+y_{\tau^{\eps,\mathbb P}_t}^\mathbb P+\eps\\
&\leq \int_t^{\tau^{\eps,\mathbb P}_t}\widehat F_s(Y_s,Z_s)ds+Y_{\tau^{\eps,\mathbb P}_t}+\eps=Y_t+ \int_t^{\tau^{\eps,\mathbb P}_t}Z_sdB_s-\int_t^{\tau_t^{\eps,\mathbb P}}{{\bf 1}_{y_{s^-}^\mathbb P<S_{s^-}}}dV_s^\mathbb P+\eps.
\end{align*}

Then, notice that the process $({\bf 1}_{y_{s^-}^\mathbb P<S_{s^-}}-{\bf 1}_{Y_{s^-} >L_{s^-}})dV_s^\mathbb P$ is non-decreasing. Therefore, we deduce
\begin{align*}
\int_t^{\sigma\wedge\tau^{\eps,\mathbb P}_t}\widehat F_s(Y_s,Z_s)ds+R_{\tau^{\eps,\mathbb P}_t}^\sigma\leq& \ Y_t+ \int_t^{\tau^{\eps,\mathbb P}_t}Z_sdB_s-\int_t^{\tau_t^{\eps,\mathbb P}}{\bf 1}_{y_{s^-}^\mathbb P<S_{s^-}}dV_s^\mathbb P\\
&+(1-\gamma)\int_t^{\tau_t^{\eps,\mathbb P}}\left({\bf 1}_{y_{s^-}^\mathbb P<S_{s^-}}-{\bf 1}_{Y_{s^-} >L_{s^-}}\right)dV_s^\mathbb P+\eps\\
=&\ Y_t+ \int_t^{\tau^{\eps,\mathbb P}_t}Z_sdB_s-(K_{\tau_t^{\mathbb P,\eps}}^\gamma-K_t^\gamma)+\eps.
\end{align*}

Similarly on the set $\{\tau_t^{\eps,\mathbb P}\geq \sigma\}$, we have
\begin{align*}
&\int_t^{\sigma\wedge\tau^{\eps,\mathbb P}_t}\widehat F_s(Y_s,Z_s)ds+R_{\tau^{\eps,\mathbb P}_t}^\sigma\leq \int_t^{\sigma}\widehat F_s(Y_s,Z_s)ds+\xi{\bf 1}_{\sigma=T}+Y_\sigma{\bf 1}_{\sigma<T}\\
&=Y_t+ \int_t^{\sigma}Z_sdB_s-\int_t^{\sigma}{\bf 1}_{y_{s^-}^\mathbb P<S_{s^-}}dV_s^\mathbb P\leq Y_t+ \int_t^{\sigma}Z_sdB_s-(K_{\sigma}^\gamma-K_t^\gamma).
\end{align*}

With these two inequalities, we therefore have
\begin{equation}
\mathbb E_t^\mathbb P\left[\int_t^{\sigma\wedge\tau^{\eps,\mathbb P}_t}\widehat F_s(Y_s,Z_s)ds+R_{\tau^{\eps,\mathbb P}_t}^\sigma+K_{\sigma\wedge\tau_t^{\mathbb P,\eps}}^\gamma-K_t^\gamma\right]-\eps\leq Y_t, \ \mathbb P-a.s.
\label{eq:dynk}
\end{equation}

We can prove similarly that for any $\tau\in\mathcal T_{t,T}$
\begin{equation}
\mathbb E_t^\mathbb P\left[\int_t^{\sigma^{\eps,\mathbb P}_t\wedge\tau}\widehat F_s(Y_s,Z_s)ds+R_{\tau}^{\sigma^{\eps,\mathbb P}_t}+K_{\sigma_t^{\mathbb P,\eps}\wedge\tau}^\gamma-K_t^\gamma\right]+\eps\geq Y_t, \ \mathbb P-a.s.
\label{eq:dynk2}
\end{equation}

Then, we can use Lemma $5.3$ of \cite{lx2} to finish the proof.
\ep

\vspace{0.35em}

Then, if we have more information on the obstacle $S$ and its decomposition \reff{eq:decomp}, we can give a more explicit representation for the processes $V^\mathbb P$, just as in the classical case (see Proposition $4.2$ in \cite{elkarquen}).

\begin{Assumption}\label{assump.s}
$S$ is a semi-martingale of the form
$$S_t=S_0+\int_0^tU_sds+\int_0^tP_sdB_s+C_t,\text{ }\mathcal P_H^\kappa-q.s.$$
where $C$ is c\`adl\`ag process of integrable variation such that the measure $dC_t$ is singular with respect to the Lebesgue measure $dt$ and which admits the following decomposition $C_t=C_t^+-C_t^-,$
where $C^+$ and $C^-$ are non-decreasing processes. Besides, $U$ and $V$ are respectively $\mathbb R$ and $\mathbb R^d$-valued $\mathcal F_t$ progressively measurable processes such that
$$\int_0^T(\abs{U_t}+\abs{P_t}^2)dt+C_T^++C_T^-\leq +\infty,\text{ }\mathcal P_H^\kappa-q.s.$$
\end{Assumption}

\vspace{0.35em}
\begin{Proposition}
Let Assumptions \ref{assump.href}, \ref{assump.obs}, \ref{assump.s} and \ref{assump.h3ref} hold. Let $(Y,Z)$ be the solution to the $2$DRBSDE \reff{2dbsderef}, then for all $\mathbb P\in\mathcal P^\kappa_H$
\begin{align}
&Z_t=P_t, \text{ }dt\times\mathbb P-a.s.\text{ on the set }\left\{Y_{t^-}=S_{t^-}\right\},
\end{align}
and there exists a progressively measurable process $(\alpha_t^\mathbb P)_{0\leq t\leq T}$ such that $0\leq \alpha\leq 1$ and
$$-{\bf 1}_{y_{t^-}^\mathbb P=S_{t^-}}dV_t^\mathbb P=\alpha_t^\mathbb P{\bf 1}_{y^\mathbb P_{t^-}=S_{t^-}}\left(\left[\widehat F_t(S_t,P_t)+U_t\right]^+dt+dC_t^+\right).$$
\end{Proposition}

\proof
First, for all $\mathbb P\in\mathcal P^\kappa_H$, the following holds $\mathbb P-a.s.$
\begin{align*}
S_t-Y_t&=S_0-Y_0+\int_0^t\left(\widehat F_s(Y_s,Z_s)+U_s\right)ds-\int_0^t(Z_s-P_s)dB_s+V_t^\mathbb P+C_t^+-C_t^-.
\end{align*}

Now if we denote $L_t$ the local time at $0$ of $S_t-Y_t$, then by It\^o-Tanaka formula under $\mathbb P$
\begin{align*}
(S_t-Y_t)^+&=(S_0-Y_0)^++\int_0^t{\bf 1}_{Y_{s^-}<S_{s^-}}\left(\widehat F_s(Y_s,Z_s)+U_s\right)ds-\int_0^t{\bf 1}_{Y_{s^-}<S_{s^-}}(Z_s-P_s)dB_s\\
&\hspace{0.9em}+\int_0^t{\bf 1}_{Y_{s^-}<S_{s^-}}d(V_t^\mathbb P+C_t^+-C_t^-)+\frac12L_t\\
&\hspace{0.9em}+\sum_{0\leq s\leq t}(S_s-Y_s)^+-(S_{s^-}-Y_{s^-})^+-{\bf 1}_{Y_{s^-}<S_{s^-}}\Delta(S_s-Y_s).
\end{align*}

However, we have $(S_t-Y_t)^+=S_t-Y_t$, hence by identification of the martingale part
$${\bf 1}_{Y_{t^-}=S_{t^-}}(Z_t-P_t)dB_t=0,\text{ }\mathcal P_H^\kappa-q.s.,$$
from which the first statement is clear. Identifying the finite variation part, we obtain
\begin{align*}
&{\bf 1}_{Y_{s^-}=S_{s^-}}\left(\widehat F_s(Y_s,Z_s)+U_s\right)ds+{\bf 1}_{Y_{s^-}=S_{s^-}}d(V_s^\mathbb P+C_s^+-C_s^-)\\
=&\frac12L_s+\left((S_s-Y_s)^+-(S_{s^-}-Y_{s^-})^+-{\bf 1}_{Y_{s^-}<S_{s^-}}\Delta(S_s-Y_s)\right).
\end{align*}

By Proposition \ref{prop.imp}, we know that ${\bf 1}_{y_{s^-}^\mathbb P=S_{s^-}}dV_s^\mathbb P$ is a non-increasing process, while ${\bf 1}_{y_{s^-}^\mathbb P<S_{s^-}}dV_s^\mathbb P$ is a non-decreasing process. Furthermore, we have
$${\bf 1}_{Y_{s^-}=S_{s^-}}dV_s^\mathbb P={\bf 1}_{y^\mathbb P_{s^-}=S_{s^-}}dV_s^\mathbb P+{\bf 1}_{y_{s^-}^\mathbb P<Y_{s^-}=S_{s^-}}dV_s^\mathbb P.$$

Since we also know that the jump part, $L$ and $C^-$ are non-decreasing processes, we obtain
\begin{align*}
-{\bf 1}_{y^\mathbb P_{s^-}=S_{s^-}}dV_s^\mathbb P &\leq {\bf 1}_{y^\mathbb P_{s^-}=S_{s^-}}\left(\left(\widehat F_s(Y_s,Z_s)+U_s\right)ds+dC_s^+\right)\\
&\hspace{0.9em}+{\bf 1}_{y_{s^-}^\mathbb P<Y_{s^-}=S_{s^-}}\left(\left(\widehat F_s(Y_s,Z_s)+U_s\right)ds+dC_s^++dV_s^\mathbb P\right).
\end{align*}

Since, ${\bf 1}_{y^\mathbb P_{s^-}=S_{s^-}}dV_s^\mathbb P$ and ${\bf 1}_{y_{s^-}^\mathbb P<Y_{s^-}=S_{s^-}}\left(\left(\widehat F_s(Y_s,Z_s)+U_s\right)ds+dC_s^++dV_s^\mathbb P\right)$ never act at the same time by definition, the second statement follows easily.
\ep

\begin{Remark}
If we assume also that $L$ is a semimartingale, then we can obtain exactly the same type of results as in Corollary 3.1 in \cite{mpz}, using the exact same arguments. 
\end{Remark}

\section{A constructive proof of existence}\label{section.3}
We have shown in Theorem \ref{uniqueref} that if a solution exists, it will necessarily verify the representation \reff{2dbsde.minref}. This gives us a natural candidate for the solution as a supremum of solutions to standard DRBSDEs. However, since those DRBSDEs are all defined on the support of mutually singular probability measures, it seems difficult to define such a supremum, because of the problems raised by the negligible sets. In order to overcome this, Soner, Touzi and Zhang proposed in \cite{stz} a pathwise construction of the solution to a 2BSDE. Let us describe briefly their strategy.

\vspace{0.35em}
The first step is to define pathwise the solution to a standard BSDE. For simplicity, let us consider first a BSDE with a generator equal to $0$. Then, we know that the solution is given by the conditional expectation of the terminal condition. In order to define this solution pathwise, we can use the so-called regular conditional probability distribution (r.p.c.d. for short) of Stroock and Varadhan \cite{str}. In the general case, the idea is similar and consists on defining BSDEs on a shifted canonical space.

\vspace{0.35em}
Finally, we have to prove measurability and regularity of the candidate solution thus obtained, and the decomposition \reff{2dbsderef} is obtained through a non-linear Doob-Meyer decomposition. Our aim in this section is to extend this approach in the presence of obstacles. We emphasize that most of the proofs are now standard, and we will therefore only sketch them, insisting particularly on the new difficulties appearing in the present setting.

\subsection{Shifted spaces}

For the convenience of the reader, we recall below some of the notations introduced in \cite{stz}.

\vspace{0.35em}
$\bullet$ For $0\leq t\leq T$, denote by $\Omega^t:=\left\{\omega\in C\left([t,T],\mathbb R^d\right),Ê\text{ }w(t)=0\right\}$ the shifted canonical space, $B^t$ the shifted canonical process, $\mathbb P_0^t$ the shifted Wiener measure and $\mathbb F^t$ the filtration generated by $B^t$.

\vspace{0.35em}
$\bullet$ For $0\leq s\leq t\leq T$ and $\omega\in \Omega^s$, define the shifted path $\omega^t\in \Omega^t$
$$\omega^t_r:=\omega_r-\omega_t,\text{ }\forall r\in [t,T].$$

\vspace{0.35em}
$\bullet$ For $0\leq s\leq t\leq T$ and $\omega\in \Omega^s$, $\widetilde \omega\in\Omega^t$ define the concatenation path $\omega\otimes_t\widetilde \omega\in\Omega^s$ by
$$(\omega\otimes_t\widetilde \omega)(r):=\omega_r1_{[s,t)}(r)+(\omega_t+\widetilde\omega_r)1_{[t,T]}(r),\text{ }\forall r\in[s,T].$$

\vspace{0.35em}
$\bullet$ For $0\leq s\leq t\leq T$ and a $\mathcal F^s_T$-measurable random variable $\xi$ on $\Omega^s$, for each $\omega \in\Omega^s$, define the shifted $\mathcal F^t_T$-measurable random variable $\xi^{t,\omega}$ on $\Omega^t$ by
$$\xi^{t,\omega}(\widetilde\omega):=\xi(\omega\otimes_t\widetilde \omega),\text{ }\forall \widetilde\omega\in\Omega^t.$$
Similarly, for an $\mathbb F^s$-progressively measurable process $X$ on $[s,T]$ and $(t,\omega)\in[s,T]\times\Omega^s$, the shifted process $\left\{X_r^{t,\omega},r\in[t,T]\right\}$ is $\mathbb F^t$-progressively measurable.

\vspace{0.35em}
$\bullet$ For a $\mathbb F$-stopping time $\tau$, the r.c.p.d. of $\mathbb P$ (denoted $\mathbb P^\omega_\tau$) is a probability measure on $\mathcal F_T$ such that
$$\mathbb E_\tau^{\mathbb P}[\xi](\omega)=\mathbb E^{\mathbb P^\omega_\tau}[\xi],\text{ for }\mathbb P-a.e.\ \omega.$$

It also induces naturally a probability measure $\mathbb P^{\tau,\omega}$ (that we also call the r.c.p.d. of $\mathbb P$) on $\mathcal F_T^{\tau(\omega)}$ which in particular satisfies that for every bounded and $\mathcal F_T$-measurable random variable $\xi$
$$\mathbb E^{\mathbb P^\omega_\tau}\left[\xi\right]= \mathbb E^{\mathbb P^{\tau,\omega}}\left[\xi^{\tau,\omega}\right].$$

\vspace{0.35em}
$\bullet$ We define similarly as in Section \ref{section.1} the set $\bar{\mathcal P}^{t}_S$, by restricting to the shifted canonical space $\Omega^t$, and its subset $\mathcal P^{t}_H$.

\vspace{0.35em}
$\bullet$ Finally, we define the "shifted" generator
$$\widehat F^{t,\omega}_s(\widetilde\omega,y,z):=F_s(\omega\otimes_t\widetilde\omega,y,z,\widehat a^t_s(\widetilde\omega)), \text{ }\forall (s,\widetilde\omega)\in[t,T]\times\Omega^t.$$

Notice that thanks to Lemma $4.1$ in \cite{stz2}, this generator coincides for $\mathbb P$-a.e.$\; \omega$ with the shifted generator as defined above, that is to say $F_s(\omega\otimes_t\widetilde\omega,y,z,\widehat a_s(\omega\otimes_t\widetilde\omega)).$ The advantage of the chosen "shifted" generator is that it inherits the uniform continuity in $\omega$ under the $\mathbb L^\infty$ norm of $F$.

\subsection{A first existence result when $\xi$ is in $\rm{UC_b}(\Omega)$ }

Let us define for all $\omega\in\Omega$, $\Lambda^*\left(\omega\right):=\underset{0\leq s\leq t}{\sup}\Lambda_s\left(\omega\right),$
where
\begin{equation*}
\Lambda_t^2\left(\omega\right):=\underset{\mathbb P\in\mathcal P^{t}_H}{\sup}\mathbb E^\mathbb P\left[\abs{\xi^{t,\omega}}^2+\int_t^T|\widehat F^{t,\omega}_s(0,0)|^2ds+\underset{t\leq s\leq T}{\sup}\abs{S^{t,\omega}_s}^2+\underset{t\leq s\leq T}{\sup}\left((L^{t,\omega}_s)^+\right)^2\right].
\end{equation*}

By Assumption \ref{assump.h3ref}, we can check directly that
\begin{equation}
\Lambda_t\left(\omega\right)<\infty \text{ for all } \left(t,\omega\right)\in\left[0,T\right]\times\Omega.
\end{equation}

To prove existence, we define the following value process $X_t$ pathwise
\begin{equation}
X_t(\omega):=\underset{\mathbb P\in\mathcal P^{t}_H}{\sup}\mathcal Y^{\mathbb P,t,\omega}_t\left(T,\xi\right), \text{ for all } \left(t,\omega\right)\in\left[0,T\right]\times\Omega,
\end{equation}
where, for any $\left(t_1,\omega\right)\in\left[0,T\right]\times\Omega,\ \mathbb P\in\mathcal P^{t_1}_H,t_2\in\left[t_1,T\right]$, and any $\mathcal F_{t_2}$-measurable $\eta\in\mathbb L^{2}\left(\mathbb P\right) $, we denote $\mathcal Y^{\mathbb P,t_1,\omega}_{t_1}\left(t_2,\eta\right):= y^{\mathbb P,t_1,\omega}_{t_1}$, where $\left(y^{\mathbb P,t_1,\omega},z^{\mathbb P,t_1,\omega},k^{\mathbb P,+,t_1,\omega},k^{\mathbb P,-,t_1,\omega}\right) $ is the solution of the following DRBSDE with upper obstacle $S^{t_1,\omega}$ and lower obstacle $L^{t_1,\omega}$ on the shifted space $\Omega^{t_1} $ under $\mathbb P$
\begin{align}\label{eq.bsdeeeerefd}
\nonumber&y^{\mathbb P,t_1,\omega}_{s}=\eta^{t_1,\omega}+\int^{t_2}_{s}\widehat{F}^{t_1,\omega}_{r}\left(y^{\mathbb P,t_1,\omega}_{r},z^{\mathbb P,t_1,\omega}_{r} \right)dr-\int^{t_2}_{s}z^{\mathbb P,t_1,\omega}_{r}dB^{t_1}_{r}-k^{\mathbb P,+,t_1,\omega}_{t_2}+k^{\mathbb P,+,t_1,\omega}_{t_1}\\
&\hspace{3.4em}+k^{\mathbb P,-,t_1,\omega}_{t_2}-k^{\mathbb P,-,t_1,\omega}_{t_1},\ \mathbb P-a.s.\\[0.5em]
&\nonumber L_t^{t_1,\omega}\leq y_t^{\mathbb P,t_1,\omega}\leq S_t^{t_1,\omega},\text{ }\mathbb P-a.s.\\
&\int_{t_1}^{t_2}\left(S_{s^-}^{t_1,\omega}-y_{s^-}^{\mathbb P,t_1,\omega}\right)dk_s^{\mathbb P,+,t_1,\omega}=\int_{t_1}^{t_2}\left(y_{s^-}^{\mathbb P,t_1,\omega}-L_{s^-}^{t_1,\omega}\right)dk_s^{\mathbb P,-,t_1,\omega}=0,\text{ }\mathbb P-a.s.
\label{skokod}
\end{align}

Notice that since we assumed that $S$ was a $\mathbb P$-semimartingale for all $\mathbb P\in\mathcal P^\kappa_H$, then for all $(t,\omega)\in[0,T]\times\Omega$, $S^{t,\omega}$ is also a $\mathbb P$-semimartingale for all $\mathbb P\in\mathcal P^{t,\kappa}_H$. Furthermore, we have the decomposition
\begin{equation}
S^{t,\omega}_s=S_t^{t,\omega}+\int_t^sP_u^{t,\omega}dB^t_u+A_s^{\mathbb P,t,\omega},\ \mathbb P-a.s.,\ \text{for all }\mathbb P\in\mathcal P^{t,\kappa}_H,
\label{eq:decompshift}
\end{equation}
where $A^{\mathbb P,t,\omega}$ is a bounded variation process under $\mathbb P$. Besides, we have by Assumption \ref{assump.h3ref}
$$\zeta^{t,\omega}_H:=\underset{\mathbb P\in\mathcal P^{t,\kappa}_H}{\sup}\left(\mathbb E^{\mathbb P}\left[\underset{t\leq s\leq T}{\esup^\mathbb P}\ \mathbb E^{H,\mathbb P}_t\left[\int_t^T\abs{(\widehat a^t_s)^{1/2}P_s^{t,\omega}}^2ds+\left(A_T^{\mathbb P,t,\omega,+}\right)^2\right]\right]\right)^{1/2}<+\infty.$$

In view of the Blumenthal zero-one law, $\mathcal Y^{\mathbb P,t,\omega}_{t}\left(T,\xi\right) $ is constant for any given $\left(t,\omega\right) $ and $\mathbb P\in\mathcal P^{t}_H $. Let us now answer the question of measurability of the process $X$
\begin{Lemma}\label{unifcontd}
Let Assumptions \ref{assump.href} and \ref{assump.h3ref} hold and consider some $\xi$ in $\rm{UC_b}(\Omega)$. Then for all $\left(t,\omega\right)\in\left[0,T\right]\times\Omega$ we have $\left|X_t\left(\omega\right)\right|\leq C(1+\zeta_H^{t,\omega}+\Lambda_t\left(\omega\right)) $. Moreover, for all $\left(t,\omega,\omega'\right)\in\left[0,T\right]\times\Omega^2$,
$\left|X_t\left(\omega\right)-X_t\left(\omega'\right)\right|\leq C\rho\left(\No{\omega-\omega'}_t\right) $. Consequently, $V_t$ is $\mathcal F_t$-measurable for every $t\in\left[0,T\right]$.
\end{Lemma}

\proof
$\rm{(i)}$ For each $\left(t,\omega\right)\in\left[0,T\right]\times\Omega $, since $S^{t,\omega}$ is a semimartingale with decomposition \reff{eq:decompshift}, we know that we have
\begin{align*}
dk_s^{\mathbb P,+,t,\omega}&\leq \left(\widehat F_s^{t,\omega}(S_s^{t,\omega},P_s^{t,\omega})\right)^+ds+dA_s^{\mathbb P,t,\omega,+}\\
&\leq C\left(\abs{\widehat F_s^{t,\omega}(0)}+\abs{S_s^{t,\omega}}+\abs{(\widehat a_s^t)^{1/2}P_s^{t,\omega}}\right)ds+dA_s^{\mathbb P,t,\omega,+}.
\end{align*}

Hence,
\begin{align}\label{label}
\nonumber\mathbb E^\mathbb P\left[\left(k_T^{\mathbb P,+,t,\omega}\right)^2\right] \leq C\left(\zeta^{t,\omega}_H+\mathbb E^\mathbb P\left[\int_t^T\abs{\widehat F_s^{t,\omega}(0)}^2 ds+\underset{t\leq s\leq T}{\sup}\abs{S_s^{t,\omega}}^2\right]\right).
\end{align}

Let now $\alpha$ be some positive constant which will be fixed later and let $\eta\in(0,1)$. By It\^o's formula we have using \reff{skokod}
\begin{align*}
&e^{\alpha t}\abs{y_t^{\mathbb P,t,\omega}}^2+\int_t^Te^{\alpha s}\abs{(\widehat a^t_s)^{1/2}z_s^{\mathbb P,t,\omega}}^2ds\leq e^{\alpha T}\abs{\xi^{t,\omega}}^2+2C\int_t^Te^{\alpha s}\abs{y_s^{\mathbb P,t,\omega}}\abs{\widehat F_s^{t,\omega}(0)}ds\\
&\hspace{0.9em}+2C\int_t^T\abs{y_s^{\mathbb P,t,\omega}}\left(\abs{y_s^{\mathbb P,t,\omega}}+\abs{(\widehat a_s^t)^{1/2}z_s^{\mathbb P,t,\omega}}\right)ds-2\int_t^Te^{\alpha s}y^{\mathbb P,t,\omega}_{s^-}z_s^{\mathbb P,t,\omega}dB^t_s\\
&\hspace{0.9em}+2\int_t^Te^{\alpha s}S_{s^-}^{t,\omega}dk_s^{\mathbb P,+,t,\omega}-2\int_t^Te^{\alpha s}L_{s^-}^{t,\omega}dk_s^{\mathbb P,-,t,\omega}-\alpha\int_t^Te^{\alpha s}\abs{y_s^{\mathbb P,t,\omega}}^2ds\\
&\leq e^{\alpha T}\abs{\xi^{t,\omega}}^2+\int_t^Te^{\alpha s}\abs{\widehat F^{t,\omega}_s(0)}^2ds-2\int_t^Te^{\alpha s}y^{\mathbb P,t,\omega}_{s^-}z_s^{\mathbb P,t,\omega}dB^t_s+\eta\int_t^Te^{\alpha s}\abs{(\widehat a^t_s)^{1/2}z_s^{\mathbb P,n}}^2ds\\[0.3em]
&\hspace{0.9em}+\left(2C+C^2+\frac{C^2}{\eta}-\alpha\right)\int_t^Te^{\alpha s}\abs{y_s^{\mathbb P,t,\omega}}^2ds+2\underset{t\leq s\leq T}{\Sup}e^{\alpha s}\abs{S_s^{t,\omega}}(k_T^{\mathbb P,+,t,\omega}-k_t^{\mathbb P,+,t,\omega})\\
&\hspace{0.9em}+2\underset{t\leq s\leq T}{\Sup}e^{\alpha s}(L_s^{t,\omega})^+(k_T^{\mathbb P,-,t,\omega}-k_t^{\mathbb P,-,t,\omega}).
\end{align*}

Now choose $\alpha$ such that $\nu:=\alpha -2C-C^2-\frac{C^2}{\eta}\geq 0$. We obtain for all $\epsilon>0$
\begin{align}
\nonumber &e^{\alpha t}\abs{y_t^{\mathbb P,t,\omega}}^2+(1-\eta)\int_t^Te^{\alpha s}\abs{(\widehat a_s^t)^{1/2}z_s^{\mathbb P,t,\omega}}^2ds\leq e^{\alpha T}\abs{\xi^{t,\omega}}^2\\
\nonumber&+\int_t^T{e^{\alpha s}{\abs{\widehat F^{t,\omega}_s(0,0)}^2}}ds+\frac1\epsilon\left(\underset{t\leq s\leq T}{\Sup}e^{\alpha s}(L_s^{t,\omega})^+\right)^2+\epsilon(k_T^{\mathbb P,-,t,\omega}-k_t^{\mathbb P,-,t\omega})^2\\
&+\left(\underset{t\leq s\leq T}{\Sup}e^{\alpha s}\abs{S_s^{t,\omega}}\right)^2+(k_T^{\mathbb P,+,t,\omega}-k_t^{\mathbb P,+,t\omega})^2-2\int_t^Te^{\alpha s}y_{s^-}^{\mathbb P,t,\omega}z_s^{\mathbb P,t,\omega}dB_s^t.
\end{align}

Taking expectation and using \reff{label} yields with $\eta$ small enough
\begin{align*}
\abs{y_t^{\mathbb P,t,\omega}}^2+\mathbb E^\mathbb P\left[\int_t^T|(\widehat a_s^t)^{\frac12}z_s^{\mathbb P,t,\omega}|^2ds\right]&\leq C\left(\Lambda_t^2(\omega)
+\left(\zeta_H^{t,\omega}\right)^2\right)+\epsilon\mathbb E^\mathbb P\left[(k_T^{\mathbb P,-,t,\omega}-k_t^{\mathbb P,-,t,\omega})^2\right].
\end{align*}

Now by definition, we also have for some constant $C_0$ independent of $\epsilon$
\begin{align*}
\mathbb E^\mathbb P\left[(k_T^{\mathbb P,-,t,\omega}-k_t^{\mathbb P,-,t,\omega})^2\right]&\leq C_0\left(\Lambda_t^2(\omega)+\left(\zeta_H^{t,\omega}\right)^2+\mathbb E^\mathbb P\left[\int_t^T\abs{y_s^{\mathbb P,t,\omega}}^2ds\right]\right)\\
&\hspace{0.9em}+C_0\mathbb E^\mathbb P\left[\int_t^T\abs{(\widehat a_s^t)^{1/2}z_s^{\mathbb P,t,\omega}}^2ds\right].
\end{align*}

Choosing $\epsilon=\frac{1}{2C_0}$, Gronwall inequality then implies $|y_t^{\mathbb P,t,\omega}|^2\leq C(1+\Lambda_t(\omega)).$ The result then follows by arbitrariness of $\mathbb P$.

\vspace{0.35em}
$\rm{(ii)}$ The proof is exactly the same as above, except that one has to use uniform continuity in $\omega$ of $\xi^{t,\omega}$, $\widehat F^{t,\omega}$, $S^{t,\omega}$ and $L^{t,\omega}$. Indeed, for each $\left(t,\omega\right)\in\left[0,T\right]\times\Omega $ and $\mathbb P\in\mathcal P^{t,\kappa}_H $, let $\alpha$ be some positive constant which will be fixed later and let $\eta\in(0,1)$. By It\^o's formula we have, since $\widehat F$ is uniformly Lipschitz
\begin{align*}
&e^{\alpha t}\abs{y_t^{\mathbb P,t,\omega}-y_t^{\mathbb P,t,\omega'}}^2+\int_t^Te^{\alpha s}\abs{(\widehat a^t_s)^{1/2}(z_s^{\mathbb P,t,\omega}-z_s^{\mathbb P,t,\omega'})}^2ds\leq e^{\alpha T}\abs{\xi^{t,\omega}-\xi^{t,\omega'}}^2\\
&\hspace{0.9em}+2C\int_t^Te^{\alpha s}\abs{y_s^{\mathbb P,t,\omega}-y_s^{\mathbb P,t,\omega'}}\left(\abs{y_s^{\mathbb P,t,\omega}-y_s^{\mathbb P,t,\omega'}}+\abs{(\widehat a_s^t)^{\frac12}(z_s^{\mathbb P,t,\omega}-z_s^{\mathbb P,t,\omega'})}\right)ds\\
&\hspace{0.9em}+2C\int_t^Te^{\alpha s}\abs{y_s^{\mathbb P,t,\omega}-y_s^{\mathbb P,t,\omega'}}\abs{\widehat F^{t,\omega}_s(y_s^{\mathbb P,t,\omega},z_s^{\mathbb P,t,\omega})-\widehat F^{t,\omega'}_s(y_s^{\mathbb P,t,\omega},z_s^{\mathbb P,t,\omega})}ds\\
&\hspace{0.9em}+2\int_t^Te^{\alpha s}(y^{\mathbb P,t,\omega}_{s^-}-y_{s^-}^{\mathbb P,t,\omega'})d\left(k_s^{\mathbb P,-,t,\omega}-k_s^{\mathbb P,-,t,\omega'}-k_s^{\mathbb P,+,t,\omega}+k_s^{\mathbb P,+,t,\omega'}\right)\\
&\hspace{0.9em}-\alpha\int_t^Te^{\alpha s}\abs{y_s^{\mathbb P,t,\omega}-y_s^{\mathbb P,t,\omega'}}^2ds-2\int_t^Te^{\alpha s}(y^{\mathbb P,t,\omega}_{s^-}-y_{s^-}^{\mathbb P,t,\omega'})(z_s^{\mathbb P,t,\omega}-z_s^{\mathbb P,t,\omega'})dB^t_s\\
&\leq e^{\alpha T}\abs{\xi^{t,\omega}-\xi^{t,\omega'}}^2+\int_t^Te^{\alpha s}\abs{\widehat F^{t,\omega}_s(y_s^{\mathbb P,t,\omega},z_s^{\mathbb P,t,\omega})-\widehat F^{t,\omega'}_s(y_s^{\mathbb P,t,\omega},z_s^{\mathbb P,t,\omega})}^2ds\\
&\hspace{0.9em}+\left(\scriptstyle2C+C^2+\frac{C^2}{\eta}-\alpha\right)\int_t^Te^{\alpha s}\abs{y_s^{\mathbb P,t,\omega}-y_s^{\mathbb P,t,\omega'}}^2ds+\eta\int_t^Te^{\alpha s}\abs{(\widehat a^t_s)^{\frac12}(z_s^{\mathbb P,t,\omega}-z_s^{\mathbb P,t,\omega'})}^2ds
\\
&\hspace{0.9em}-2\int_t^Te^{\alpha s}(y^{\mathbb P,t,\omega}_{s^-}-y_{s^-}^{\mathbb P,t,\omega'})(z_s^{\mathbb P,t,\omega}-z_s^{\mathbb P,t,\omega'})dB^t_s\\
&\hspace{0.9em}+2\int_t^Te^{\alpha s}(y^{\mathbb P,t,\omega}_{s^-}-y_{s^-}^{\mathbb P,t,\omega'})d\left(k_s^{\mathbb P,-,t,\omega}-k_s^{\mathbb P,-,t,\omega'}-k_s^{\mathbb P,+,t,\omega}+k_s^{\mathbb P,+,t,\omega'}\right).
\end{align*}

By the Skorohod condition \reff{skokod}, we also have
\begin{align*}
&\int_t^Te^{\alpha s}(y^{\mathbb P,t,\omega}_{s^-}-y_{s^-}^{\mathbb P,t,\omega'})d\left(k_s^{\mathbb P,-,t,\omega}-k_s^{\mathbb P,-,t,\omega'}-k_s^{\mathbb P,+,t,\omega}+k_s^{\mathbb P,+,t,\omega'}\right)\\
&\leq \int_t^Te^{\alpha s}(L^{t,\omega}_{s^-}-L_{s^-}^{t,\omega'})d(k_s^{\mathbb P,-,t,\omega}-k_s^{\mathbb P,-,t,\omega'})-\int_t^Te^{\alpha s}(S^{t,\omega}_{s^-}-S_{s^-}^{t,\omega'})d(k_s^{\mathbb P,+,t,\omega}-k_s^{\mathbb P,+,t,\omega'}).
\end{align*}

Now choose $\alpha$ such that $\nu:=\alpha -2C-C^2-\frac{C^2}{\eta}\geq 0$. We obtain for all $\epsilon>0$
\begin{align}
&\nonumber e^{\alpha t}\abs{y_t^{\mathbb P,t,\omega}-y_t^{\mathbb P,t,\omega'}}^2+(1-\eta)\int_t^Te^{\alpha s}\abs{(\widehat a^t_s)^{1/2}(z_s^{\mathbb P,t,\omega}-z_s^{\mathbb P,t,\omega'})}^2ds\\
&\nonumber\leq  e^{\alpha T}\abs{\xi^{t,\omega}-\xi^{t,\omega'}}^2+\int_t^Te^{\alpha s}\abs{\widehat F^{t,\omega}_s(y_s^{\mathbb P,t,\omega},z_s^{\mathbb P,t,\omega})-\widehat F^{t,\omega'}_s(y_s^{\mathbb P,t,\omega},z_s^{\mathbb P,t,\omega})}^2ds\\
\nonumber&\hspace{0.9em}+\frac1\epsilon\left(\underset{t\leq s\leq T}{\Sup}e^{\alpha s}(L_s^{t,\omega}-L_s^{t,\omega'})^+\right)^2+\epsilon(k_T^{\mathbb P,-,t,\omega}-k_T^{\mathbb P,-,t,\omega'}-k_t^{\mathbb P,-,t\omega}+k_t^{\mathbb P,-,t,\omega'})^2\\
\nonumber&\hspace{0.9em}+\frac1\epsilon\left(\underset{t\leq s\leq T}{\Sup}e^{\alpha s}\abs{S_s^{t,\omega}-S_s^{t,\omega'}}\right)^2+\epsilon(k_T^{\mathbb P,+,t,\omega}-k_T^{\mathbb P,+,t,\omega'}-k_t^{\mathbb P,+,t\omega}+k_t^{\mathbb P,+,t,\omega'})^2\\
&\hspace{0.9em}-2\int_t^Te^{\alpha s}(y^{\mathbb P,t,\omega}_{s^-}-y_{s^-}^{\mathbb P,t,\omega'})(z_s^{\mathbb P,t,\omega}-z_s^{\mathbb P,t,\omega'})dB^t_s.
\label{bidule3}
\end{align}

The end of the proof is then similar to the previous step, using the uniform continuity in $\omega$ of $\xi$, $F$ and $S$.
\ep

\vspace{0.35em}
Then, we show the same dynamic programming principle as Proposition $4.7$ in \cite{stz2} and Proposition $4.1$ in \cite{mpz}. The proof being exactly the same, we omit it.
\begin{Proposition}\label{progdynd}
Under Assumptions \ref{assump.href}, \ref{assump.h3ref} and for $\xi\in \rm{UC_b}(\Omega)$, we have for all $0\leq t_1<t_2\leq T$ and for all $\omega \in \Omega$
$$X_{t_1}(\omega)=\underset{\mathbb P\in \mathcal P^{t_1,\kappa}_H}{\sup}\mathcal Y_{t_1}^{\mathbb P,t_1,\omega}(t_2,X_{t_2}^{t_1,\omega}).$$
\end{Proposition}

\vspace{0.35em}
Define now for all $(t,\omega)$, the $\mathbb F^+$-progressively measurable process
\begin{equation}
\label{vplus}X_t^+:=\underset{r\in\mathbb Q\cap(t,T],r\downarrow t}{\overline \lim}X_r.
\end{equation} We have the following result whose proof is the same as the one of Lemma $4.2$ in \cite{mpz}.

\begin{Lemma}
Under the conditions of the previous Proposition, we have
$$X_t^+=\underset{r\in\mathbb Q\cap(t,T],r\downarrow t}{\lim}X_r,\text{ }\mathcal P^\kappa_H-q.s.$$
and thus $X^+$ is c\`adl\`ag $\mathcal P^\kappa_H-q.s.$
\end{Lemma}

\vspace{0.35em}
Proceeding exactly as in Steps $1$ et $2$ of the proof of Theorem $4.5$ in \cite{stz2}, we can then prove that $X^+$ is a strong doubly reflected $\widehat F$-supermartingale (in the sense of Definition \ref{defdoob} in the Appendix). Then, using the Doob-Meyer decomposition proved in the Appendix in Theorem \ref{doobmeyerd} for all $\mathbb P$, we know that there exists a unique ($\mathbb P-a.s.$) process $\overline Z^\mathbb P\in\mathbb H^2(\mathbb P)$ and unique non-decreasing c\`adl\`ag square integrable processes $A^\mathbb P$, $B^\mathbb P$ and $C^\mathbb P$ such that

\begin{itemize}
\item $X_t^+=X_0^+-\int_0^t\widehat F_s(X_s^+,\overline Z_s^\mathbb P)ds+\int_0^t\overline Z^\mathbb P_sdB_s+B_t^\mathbb P-A_t^\mathbb P+C_t^\mathbb P, \text{ }\mathbb P-a.s., \text{ }\forall \mathbb P\in\mathcal P^\kappa_H.$
\item $L_t\leq X_t^+\leq S_t, \text{ }\mathbb P-a.s., \text{ }\forall \mathbb P\in\mathcal P^\kappa_H.$
\item $\int_0^T\left(S_{t^-}-X_{t^-}^+\right)dA_t^\mathbb P=\int_0^T\left(X_{t^-}^+-L_{t^-}\right)dB_t^\mathbb P=0, \text{ }\mathbb P-a.s., \text{ }\forall \mathbb P\in\mathcal P^\kappa_H.$
\end{itemize}

\vspace{0.35em}
We then define $V^\mathbb P:=A^\mathbb P-B^\mathbb P-C^\mathbb P$. By Karandikar \cite{kar}, since $X^+$ is a c\`adl\`ag semimartingale, we can define a universal process $\overline Z$ which aggregates the family $\left\{\overline Z^\mathbb P,\mathbb P\in\mathcal P^\kappa_H\right\}$.

\vspace{0.35em}
We next prove the representation \reff{representationref} for $X$ and $X^+$.
\begin{Proposition}\label{prop.reprefd}
Assume that $\xi\in UC_b(\Omega)$ and that Assumptions \ref{assump.href} and \ref{assump.h3ref} hold. Then we have
$$X_t=\underset{\mathbb P^{'}\in\mathcal P_H^\kappa(t,\mathbb P)}{\esup^\mathbb P}\mathcal Y_t^{\mathbb P^{'}}(T,\xi)\text{ and } X_t^+=\underset{\mathbb P^{'}\in\mathcal P^\kappa_H(t^+,\mathbb P)}{\esup^\mathbb P}\mathcal Y_t^{\mathbb P^{'}}(T,\xi), \text{ }\mathbb P-a.s., \text{ }\forall \mathbb P\in\mathcal P_H^\kappa.$$
\end{Proposition}

\proof
The proof for the representations is the same as the proof of proposition $4.10$ in \cite{stz2}, since we also have a stability result for RBSDEs under our assumptions.
\ep

\vspace{0.35em}
Finally, we have to check that the minimum condition \reff{2dbsde.minref} holds. However, this can be done exactly as in \cite{mpz}, so we refer the reader to that proof.

\subsection{Main result}
We are now in position to state the main result of this section

\begin{Theorem}\label{mainrefd}
Let $\xi\in\mathcal L^{2,\kappa}_H$ and let Assumptions \ref{assump.href},\ref{assump.obs} and \ref{assump.h3ref} hold. Then

\vspace{0.35em}
1)   There exists a unique solution $(Y,Z)\in\mathbb D^{2,\kappa}_H\times\mathbb H^{2,\kappa}_H$ of the $2\rm{DRBSDE}$ \reff{2dbsderef}. \\
2) Moreover, if in addition we choose to work under either of the following model of set theory (we refer the reader to \cite{frem} for more details)
\begin{itemize}
\item[\rm{(i)}] Zermelo-Fraenkel set theory with axiom of choice (ZFC) plus the Continuum Hypothesis (CH).
\item[\rm{(ii)}] ZFC plus the negation of CH plus Martin's axiom.
\end{itemize}
Then there exists a unique solution $(Y,Z,V)\in\mathbb D^{2,\kappa}_H\times\mathbb H^{2,\kappa}_H\times\mathbb V^{2,\kappa}_H$ of the $2\rm{DRBSDE}$ \reff{2dbsderef}.
\end{Theorem}

\proof
The proof of the existence part follows the lines of the proof of Theorem $4.7$ in \cite{stz}, using the estimates of Theorem \ref{estimates2d}, so we only insist on the points which do not come directly from the proofs mentioned above. The idea is to approximate the terminal condition $\xi$ by a sequence $(\xi^n)_{n\geq 0}\subset {\rm{UC}_b(\Omega)}$. Then, we use the estimates of Theorem \ref{estimates2d} to pass to the limit as in the proof of Theorem $4.6$ in \cite{stz}. The main point in this context is that for each $n$, if we consider the Jordan decomposition of $V^{\mathbb P,n}$ into the non-decreasing process $V^{\mathbb P,+,n}$ and the non-increasing process $V^{\mathbb P,-,n}$, then the estimates of Theorem \ref{estimates2d} ensure that these processes converge to some $V^{\mathbb P,+}$ and $V^{\mathbb P,-}$, which are respectively non-decreasing and non-increasing. Hence we are sure that the limit $V^\mathbb P$ has indeed bounded variation. 

\vspace{0.5em}
Concerning the fact that we can aggregate the family $\left(V^\mathbb P\right)_{\mathbb P\in\mathcal P^\kappa_H}$, it can be deduced as follows.
\vspace{0.5em}
First, if $\xi\in \rm{UC}_b(\Omega)$, we know, using the same notations as above that the solution verifies
$$X_t^+=X_0^+-\int_0^t\widehat F_s(X_s^+,\overline Z_s)ds+\int_0^t\overline Z_sdB_s-K_t^\mathbb P, \text{ }\mathbb P-a.s., \text{ }\forall \mathbb P\in\mathcal P_H^\kappa.$$

Now, we know from \reff{vplus} that $X^+$ is defined pathwise, and so is the Lebesgue integral $$\int_0^t\widehat F_s(X_s^+,\overline Z_s)ds.$$

In order to give a pathwise definition of the stochastic integral, we would like to use the recent results of Nutz \cite{nutz}. However, the proof in this paper relies on the notion of medial limits, which may or may not exist depending on the model of set theory chosen. They exists in the model (i) above, which is the one considered by Nutz, but we know from \cite{frem} (see statement $22$O(l) page $55$) that they also do in the model (ii). Therefore, provided we work under either one of these models, the stochastic integral $\int_0^t\overline Z_sdB_s$ can also be defined pathwise. We can therefore define pathwise
$$V_t:=X_0^+-X_t^+-\int_0^t\widehat F_s(X_s^+,\overline Z_s)ds+\int_0^t\overline Z_sdB_s,$$
and $V$ is an aggregator for the family $\left(V^\mathbb P\right)_{\mathbb P\in\mathcal P^\kappa_H}$, that is to say that it coincides $\mathbb P-a.s.$ with $V^\mathbb P$, for every $\mathbb P\in\mathcal P^\kappa_H$.

\vspace{0.5em}
In the general case when $\xi\in\mathcal L^{2,\kappa}_H$, the family is still aggregated when we pass to the limit.
\ep


\begin{Remark}
For more discussions on the axioms of the set theory considered here, we refer the reader to the Remark $4.2$ in \cite{mpz}.\end{Remark}

\section{Applications: Israeli options and Dynkin games}\label{sec.app}

\subsection{Game options}

We first recall the definition of an Israeli (or game) option, and we refer the reader to \cite{kif}, \cite{ham2} and the references therein for more details. An Israeli option is a contract between a broker (seller) and a trader (buyer). The specificity is that both can decide to exercise before the maturity date $T$. If the trader exercises first at a
time $t$ then the broker pays him the (random) amount $L_t$. If the
broker exercises before the trader at time $t$, the trader will be
given from him the quantity $S_t\geq L_t$, and the difference $S_t-L_t$ as to be understood as a penalty imposed on the seller for canceling the contract. In the case where
they exercise simultaneously at $t$, the trader payoff is $L_t$
and if they both wait till the maturity of the contract $T$, the trader receives the amount $\xi$. In other words, this is an American option which has the specificity that the seller can also "exercise" early. This therefore a typical Dynkin game. We assume throughout this section that the processes $L$ and $S$ satisfy Assumptions \ref{assump.s} and \ref{assump.h3ref}.

\vspace{0.35em}
To sum everything up, if we consider that the broker exercises at
a stopping time $\tau\leq T$ and the trader at another time $\sigma\leq T$
then the trader receive from the broker the following payoff:
$$H(\sigma,\tau):=S_\tau 1_{\tau<\sigma}+L_\sigma 1_{\sigma\leq \tau}+\xi 1_{\sigma\wedge\tau=T}$$

\begin{Remark} We could have chosen a slightly more general payoff function as in \cite{ham2}, but we prefer to concentrate here on the uncertainty context.
\end{Remark}

Before introducing volatility uncertainty, let us first briefly recall how the fair price and the hedging of such an option is related to DRBSDEs in a classical financial market. We fix a probability measure $\mathbb P$, and we assume that the market contains one riskless asset, whose price is assumed w.l.o.g. to be equal to $1$, and one risky asset. We furthermore assume that if the broker adopts a strategy $\pi$ (which is an adapted process in $\mathbb H^2(\mathbb P)$ representing the percentage of his total wealth invested in the risky asset), then his wealth process has the following expression
 $$X_t^\P=\xi+\int_t^Tb(s,X_s^\P,\pi_s^\P)ds-\int_t^T\pi_s^\P\sigma_s dW_s\, ,\P-a.s.$$

where W is a Brownian motion under $\P$, $b$ is convex and Lipschitz
with respect to $(x,\pi)$. We also suppose that the process
$(b(t,0,0))_{t\leq T}$ is square-integrable and $(\sigma_t)_{t\leq T}$ is invertible and its inverse is bounded. It was then proved in \cite{kif} and \cite{ham2} that the fair price and an hedging strategy for the Israeli option described above can be obtained through the solution of a DRBSDE. More precisely, we have

\begin{Theorem} The fair price of the {\it game option} and the corresponding hedging strategy are given by the pair $(y^\mathbb P,\pi^\P)\in\mathbb D^2(\mathbb P)\times\H^2(\P)$ solving the following DRBSDE
$$\begin{cases}
\nonumber y_t^\mathbb P=\xi+\int_t^Tb(s,y_s^\P,\pi_s^\P)ds-\int_t^T\pi_s^\P\sigma_s dW_s +k^\P_t-k^\P_t,\, \mathbb P-a.s.\\
\nonumber L_t\leq y_t^\P\leq S_t, \, \mathbb P-a.s.\\
\int_0^T(y_{t^-}^\P-L_{t^-})dk^{\P,-}_t=\int_0^T(S_{t^-}-y_{t^-}^\P)dk^{\P,+}_t=0.
\end{cases}$$
Moreover, for any $\eps>0$, the following stopping times are $\eps$-optimal after $t$ for the seller and the buyer respectively
$$D^{1,\eps,\mathbb P}_t:=\inf\left\{s\geq t,\ y_s^\mathbb P\geq S_s-\eps\right\}, \ D^{2,\eps,\mathbb P}_t:=\inf\left\{s\geq t,\ y_s^\mathbb P\leq L_s+\eps\right\}.$$
\end{Theorem}

Let us now extend this result to the uncertain volatility framework.
We still consider a financial market with two assets and assume now
that the wealth process has the following dynamic when the chosen strategy is $\pi$
$$X_t=\xi+\int_t^Tb(s,X_s,\pi_s)ds-\int_t^T\pi_s dB_s,\,\mathcal{P}^\kappa_H-q.s., $$
where $B$ is the canonical process and $b$ is assumed to satisfy Assumptions \ref{assump.href}, and where $\xi$ belongs to $\mathcal{L}^{2,\kappa}_H$. Then, following the ideas of \cite{mpz}, it is natural to consider as a superhedging price for the option the quantity
$$ Y_t=\underset{\P'\in\mathcal{P}_H^\kappa(t^+,\P)}{\esup^\mathbb P} y_t^{\P'}.$$

Indeed, this amount is greater than the price at time $t$ of the the same Israeli option under any probability measure. Hence, if the seller receives this amount, he should always be able to hedge his position. We emphasize however that we are not able to guarantee that this price is optimal in the sense that it is the lowest value for
which we can find a super-replicating strategy. This interesting question is left for future research.

\vspace{0.35em}
Symmetrically, if the seller charges less than the following quantity for the option at time $t$,
$$\widetilde Y_t:=\underset{\P'\in\mathcal{P}_H^\kappa(t^+,\P)}\einf^\mathbb P y_t^{\P'},$$
then it will be clearly impossible for him to find a hedge. $\widetilde Y$ appears then as a subhedging price.

\vspace{0.35em}
Hence, we have obtained a whole interval of prices, given by $[\widetilde Y_t,Y_t]$, which we can formally think as arbitrage free, even though a precise definition of this notion in an uncertain market is outside the scope of this paper. These two quantities can be linked to the notion of second-order 2DRBSDEs. Indeed, this is immediate for $Y$, and for $\widetilde Y$, we need to introduce a "symmetric" definition for the 2DRBSDEs.

\begin{Definition}\label{death3}
For $\xi \in \mathbb L^{2,\kappa}_H$, we consider the following type of equations satisfied by a pair of progressively-measurable processes $(Y,Z)$
\begin{itemize}
\item[$\bullet$] $Y_T=\xi$, $\mathcal P_H^\kappa-q.s$.
\item[$\bullet$] $\forall \mathbb P \in \mathcal P_H^\kappa$, the process $V^{\mathbb P}$ defined below has paths of bounded variation $\mathbb P-a.s.$
\begin{equation}
V_t^{\mathbb P}:=Y_0-Y_t - \int_0^t\widehat{F}_s(Y_s,Z_s)ds+\int_0^tZ_sdB_s, \text{ } 0\leq t\leq T, \text{  } \mathbb P-a.s.
\label{2dbsde.kref2}
\end{equation}

\item[$\bullet$] We have the following maximum condition for $0\leq t\leq T$
\begin{equation}
V_t^\mathbb P+k_t^{\mathbb P,+}-k_t^{\mathbb P,-}=\underset{ \mathbb{P}^{'} \in \mathcal{P}_H(t^+,\mathbb{P}) }{ \esup^{\mathbb P} }\mathbb{E}_t^{\mathbb P^{'}}\left[V_T^{\mathbb P^{'}}+k_T^{\mathbb P^{'},+}-k_T^{\mathbb P^{'},-}\right], \text{  } \mathbb P-a.s., \text{ } \forall \mathbb P \in \mathcal P_H^\kappa.
\label{2dbsde.minref2}
\end{equation}
\item[$\bullet$] $L_t\leq Y_t\leq S_t$, $\mathcal P_H^\kappa-q.s.$
\end{itemize}
\end{Definition}

This Definition is symmetric to Definition \ref{death2} in the sense that if $(Y,Z)$ solves an equation as in Definition \ref{death3}, then $(-Y,-Z)$ solves a 2DRBSDE (in the sense of Definition \ref{death2}) with terminal condition $-\xi$, generator $\tilde g(y,z):=-g(-y,-z)$, lower obstacle $-S$ and upper obstacle $-L$. With this remark, it is clear that we can deduce a wellposedness theory for the above equations. In particular, we have the following representation
\begin{equation}\label{rep2}
Y_t=\underset{\mathbb P^{'}\in\mathcal P^\kappa_H(t^+,\mathbb P)}{\einf}y_t^{\mathbb P^{'}},\ \P-a.s.,\text{ for any }\mathbb P\in\mathcal P^\kappa_H.
\end{equation}

We then have the following result.

\begin{Theorem} The superhedging and subhedging prices $Y$ and $\widetilde Y$ are respectively the unique solution of the 2DRBSDE with terminal condition $\xi$, generator $b$, lower obstacle $L$, upper obstacle $S$ in the sense of Definitions \ref{death2} and \ref{death3} respectively. The corresponding hedging strategies are then given by $Z$ and $\widetilde Z$.

\vspace{0.35em}
Moreover, for any $\eps>0$ and for any $\mathbb P$, the following stopping times are $\eps$-optimal after $t$ for the seller and the buyer respectively
$$D^{1,\eps,\mathbb P}_t:=\inf\left\{s\geq t,\ y_s^\mathbb P\geq S_s-\eps,\ \mathbb P-a.s.\right\}, \ D^{2,\eps}_t:=\inf\left\{s\geq t,\ Y_s\leq L_s+\eps,\ \mathcal P^\kappa_H-q.s.\right\}.$$
\end{Theorem}

\subsection{A first step towards Dynkin games under uncertainty}
It is already known that doubly reflected BSDE are intimately connected to
Dynkin games (see \cite{cvi} for instance). More generally, since the seminal paper by Fleming and Souganidis \cite{flem}, two person zero-sum stochastic differential games have been typically studied through two approaches. One uses the viscosity theory and aims at showing that the value function of the game is the unique viscosity solution of the associated HJBI equation, while the other relates the value function to the solution of a BSDE. We are of course more interested in the second one. To name but a few of the contributions in the literature, Buckdahn and Li \cite{buck} defined precisely the value function of the game via BSDEs, while more recently Bayraktar and Yao used doubly reflected BSDEs. Before specializing the discussion to Dynkin games, we would like to refer the reader to the very recent work of Pham and Zhang \cite{pham}, which studies a weak formulation of two person zero-sum game and points out several formal connections with the 2BSDE theory.

\vspace{0.35em}
We naturally want to obtain the same kind of result with 2DRBSDEs, with an additional uncertainty component in the
game, which will be induced by the fact that we are working simultaneously under a family of mutually singular probability measures. We will focus here
on the construction of a game whose
upper and lower values can be expressed as a solution of 2DRBSDE. We insist that we prove only that a given solution to a 2DRBSDE provides a solution to the corresponding Dynkin game described below. However, we are not able, as in \cite{cvi}, to construct the solution of the 2DRBSDE directly from the solution of the Dynkin game. Moreover, we also face a difficult technical problem related to Assumption \ref{assump.minmax}, which prevents our result to be comprehensive.

\vspace{0.35em}
Let us now describe what we mean precisely by a Dynkin game with uncertainty. Two players P1 and P2 are facing each other in a game. A strategy of
a player consists in picking a stopping time. Let us say that P1 chooses
$\tau\in\mathcal T_{0,T}$ and P2 chooses $\sigma\in\mathcal T_{0,T}$. Then the game stipulates that P1 will pay to
P2 the following random payoff
$$R_t(\tau,\sigma):=
\int_t^{\tau\wedge\sigma}g_sds+S_\tau1_{\tau<\sigma}+L_\sigma1_{\sigma\leq\tau,\sigma<T}+\xi1_{\tau\wedge\sigma=T},$$
where $g$, $S$ and $L$ are $\mathbb F$-progressively measurable processes satisfying Assumptions \ref{assump.href}, \ref{assump.s} and \ref{assump.h3ref}. In particular, the upper obstacle $S$ is a semimartingale.

\vspace{0.35em}
Then naturally, P1 will try to minimize the expected amount that he will have to pay, but taking into account the fact that both P2 and the "Nature" (which we interpret as a third player, represented by the uncertainty that the player has with respect to the underlying probability measure) can play against him. Symmetrically, P2 will try to maximize his expected returns, considering that both P1 and the Nature are antagonist players. This leads us to introduce the following upper and lower values of the robust Dynkin game
$$\overline V_t:=\underset{\tau\in\mathcal{T}_{t,T}}\einf\, \underset{\sigma\in\mathcal{T}_{t,T}}\esup\,\underset{\P'\in\mathcal{P}^\kappa_H(t^+,\P)}{\esup^\P}\, \E_t^{\P'}[R_t(\tau,\sigma)],\ \mathbb P-a.s.$$
$$\underline{V}_t:=\underset{\sigma\in\mathcal{T}_{t,T}}\esup\, \underset{\tau\in\mathcal{T}_{t,T}}\einf\,\underset{\P'\in\mathcal{P}^\kappa_H(t^+,\P)}{\einf^\P}\, \E_t^{\P'}[R_t(\tau,\sigma)],\ \mathbb P-a.s.$$

\begin{Remark}
In order to be completely rigorous, we should have made the dependence in $\mathbb P$ of the two functions above explicit, because it is not clear that an aggregator exists a priori. Nonetheless, we will prove in this section, that they both correspond to the solution of 2DRBSDE, and therefore that the aggregator indeed exists. Therefore, for the sake of clarity, we will always omit this dependence.
\end{Remark}

$\overline V$ is the maximal amount that P1 will agree to pay in order to take part in the game. Symmetrically, $\underline V$ is the minimal amount that P2 must receive in order to accept to take part to the game. Unlike in the classical setting without uncertainty, for which there is only one value on which the 2 players can agree, in our context there is generally a whole interval of admissible values for the game. Indeed, we have the following easy result

\begin{Lemma}
We have for $t\in [0,T]$
$$\overline V_t\geq \underline V_t,\ \mathcal P^\kappa_H-q.s.$$
Therefore the admissible values for the game are the interval $[\underline V_t,\overline V_t]$.
\end{Lemma}

\proof
Let $t\in[0,T]$, $\mathbb P\in\mathcal P^\kappa_H$ and $\mathbb P^{'}\in\mathcal P^\kappa_H(t^+,\mathbb P)$. For any $(\tau,\sigma)\in\mathcal T_{t,T}\times\mathcal T_{t,T}$, we have clearly
$$\underset{\P'\in\mathcal{P}^\kappa_H(t^+,\P)}{\esup^\P}\, \E_t^{\P'}[R_t(\tau,\sigma)]\geq \underset{\tau\in\mathcal{T}_{t,T}}\einf\,\underset{\P'\in\mathcal{P}^\kappa_H(t^+,\P)}{\einf^\P}\, \E_t^{\P'}[R_t(\tau,\sigma)],\ \mathbb P-a.s.$$

Then we can take the essential supremum with respect to $\sigma$ on both sides of the inequality, and the result follows.
\ep

\vspace{0.35em}
Now, in order to link the solution of the above robust Dynkin game to 2DRBSDEs, we will need to assume a min-max property which is closely related to the usual Isaacs condition for the classical Dynkin games. Given the length of the paper, we will not try to verify this assumption. Nonetheless, we emphasize that a related result was indeed proved in \cite{mpz} in the context of a robust utility maximization problem. Even more in the spirit of our paper, Nutz and Zhang \cite{nutz2} also showed such a result (at least at time $t=0$ and under sufficient regularity assumptions, see their Theorem $3.4$) when there is only one player. We are convinced that their results could be generalized to our framework and leave this interesting problem to future research.

\begin{Assumption}\label{assump.minmax}
We suppose that the following "min-max" property are satisfied. For any $\mathbb P\in\mathcal P^\kappa_H$
\begin{align}\label{dyn1}
&\underset{\tau\in\mathcal{T}_{t,T}}\einf\, \underset{\sigma\in\mathcal{T}_{t,T}}\esup\,\underset{\P'\in\mathcal{P}^\kappa_H(t^+,\P)}{\esup^\P}\, \E_t^{\P'}[R_t(\tau,\sigma)]
=\underset{\P'\in\mathcal{P}^\kappa_H(t^+,\P)}{\esup^\P}\,\underset{\tau\in\mathcal{T}_{t,T}}\einf\,
\underset{\sigma\in\mathcal{T}_{t,T}}\esup\,
\E_t^{\P'}[R_t(\tau,\sigma)],\ \mathbb P-a.s.\\
&\underset{\sigma\in\mathcal{T}_{t,T}}\esup\, \underset{\tau\in\mathcal{T}_{t,T}}\einf\,\underset{\P'\in\mathcal{P}^\kappa_H(t^+,\P)}{\einf^\P}\, \E_t^{\P'}[R_t(\tau,\sigma)]
=\underset{\P'\in\mathcal{P}^\kappa_H(t^+,\P)}{\einf^\P}\,\underset{\sigma\in\mathcal{T}_{t,T}}\esup\,
\underset{\tau\in\mathcal{T}_{t,T}}\einf\,
\E_t^{\P'}[R_t(\tau,\sigma)],\ \mathbb P-a.s.
\label{dyn2}
\end{align}
\end{Assumption}

It is clear from Proposition \ref{dynkin} that the right-hand side of \reff{dyn1} can be expressed as the solution of 2DRBSDE with terminal condition $\xi$, generator $g$, lower obstacle $L$ and upper obstacle $S$.
We deduce immediately the following result.

\begin{Theorem}
Let Assumption \ref{assump.minmax} hold. Let $(Y,Z)$ $($resp. $(\widetilde Y,\widetilde Z)$$)$ be a solution to the 2DRBSDE in the sense of Definition \ref{death2} (resp. in the sense of Definition \ref{death3}) with terminal condition $\xi$, generator $g$, lower obstacle $L$ and upper obstacle $S$. Then we have for any $t\in[0,T]$
\begin{align*}
&\overline V_t=Y_t,\ \mathcal P^\kappa_H-q.s.\\
&\underline V_t=\widetilde Y_t,\ \mathcal P^\kappa_H-q.s.
\end{align*}
Moreover, unless $\mathcal P^\kappa_H$ is reduced to a singleton, we have $\overline{V}>\underline{V},\ \mathcal P^\kappa_H-q.s.$
\end{Theorem}

\proof
The two equalities are obvious. Moreover, if for each $\mathbb P\in\mathcal P^\kappa_H$, we let $y^\mathbb P$ be the solution of the DRBSDE with terminal condition $\xi$, generator $g$, lower obstacle $L$ and upper obstacle $S$, we have by \reff{rep2}
$$\overline V_t=\underset{\mathbb P^{'}\in\mathcal P^\kappa_H(t^+,\mathbb P)}{\esup}y_t^{\mathbb P^{'}},\ \P-a.s., \text{ and }\underline V_t=\underset{\mathbb P^{'}\in\mathcal P^\kappa_H(t^+,\mathbb P)}{\einf}y_t^{\mathbb P^{'}},\ \P-a.s.,$$
which implies the last result.
\ep

\begin{appendix}

\section{Doubly reflected $g$-supersolution and martingales}\label{sec.sec}
In this section, we extend some of the results of \cite{pengg} and \cite{mpz} concerning $g$-supersolution of BSDEs and RBSDEs to the case of DRBSDEs. Let us note that many of the results below are obtained using similar ideas as in \cite{pengg} and \cite{mpz}, but we still provide most of them since, to the best of our knowledge, they do not appear anywhere else in the literature. Moreover, we emphasize that we only provide the results and definitions for the doubly reflected case, because the corresponding ones for the upper reflected case can be deduced easily. In the following, we fix a probability measure $\mathbb P$

\subsection{Definitions and first properties}\label{sec.sec.sec}
Let us be given the following objects

\vspace{0.35em}
$\bullet$ A function $g_s(\omega,y,z)$, $\mathbb F$-progressively measurable for fixed $y$ and $z$, uniformly Lipschitz in $(y,z)$ and such that $\mathbb E^\mathbb P\left[\int_0^T\abs{g_s(0,0)}^2ds\right]<+\infty.$

\vspace{0.35em}
$\bullet$ A terminal condition $\xi$ which is $\mathcal F_T$-measurable and in $L^2(\mathbb P)$.

\vspace{0.35em}
$\bullet$ C\`adl\`ag processes $V$, $S$, $L$ in $\mathbb I^2(\mathbb P)$ such that $S$ and $L$ satisfy Assumptions \ref{assump.obs} and \ref{assump.h3ref} (in particular $S$ is a semimartingale with the decomposition \reff{eq:decomp}), and with $\mathbb E^\mathbb P\left[\underset{0\leq t\leq T}{\sup}\abs{V_t}^2\right]<+\infty$.

\vspace{0.55em}
We study the problem of finding $(y,z,k^+,k^-)\in\mathbb D^2(\mathbb P)\times\mathbb H^2(\mathbb P)\times\mathbb I^2(\mathbb P)\times\mathbb I^2(\mathbb P)$ such that
\begin{align}
\label{drbsde2}
\nonumber&y_t=\xi+\int_t^Tg_s(y_s,z_s)ds-\int_t^T z_sdW_s+k^-_T-k^-_t-k^+_{T}+k^+_{t}+V_T-V_t,\text{ } \mathbb P-a.s.\\
\nonumber&L_t\leq y_t\leq S_t, \text{ }\mathbb P-a.s.\\
&\int_0^T\left(S_{s^-}-y_{s^-}\right)dk^+_s=\int_0^T\left(y_{s^-}-L_{s^-}\right)dk^-_s=0, \text{ }\mathbb P-a.s.
\end{align}

We first have a result of existence and uniqueness.

\begin{Proposition}\label{ex}
Under the above hypotheses, there is a unique solution $(y,z,k^+,k^-)\in\mathbb D^2(\mathbb P)\times\mathbb H^2(\mathbb P)\times\mathbb I^2(\mathbb P)\times\mathbb I^2(\mathbb P)$ to the doubly reflected BSDE \reff{drbsde2}.
\end{Proposition}

\proof
Consider the following penalized RBSDE with lower obstacle $L$, whose existence and uniqueness are ensured by the results of Lepeltier and Xu \cite{lx}
$$y^{n}_t=\xi+\int_t^Tg_s(y^{n}_s,z^{n}_s)ds-\int_t^Tz^{n}_sdW_s+k_T^{n,-}-k_t^{n,-}-k_T^{n,+}+k_t^{n,+}+V_T-V_t,$$
where $k_t^{n,+}:=n\int_0^t(S_s-y^{n}_s)^-ds$. Then, define $\widetilde y_t^{n}:=y_t^{n}+V_t$, $\widetilde \xi:=\xi+V_T$, $\widetilde z_t^{n}:=z_t^{n}$, $\widetilde k_t^{n,\pm}:=k_t^{n,\pm}$, $\widetilde g_t(y,z):=g_t(y-V,z)$ and $\widetilde L_t:=L_t+V_t$. Then
$$\widetilde y^{n}_t=\widetilde \xi+\int_t^T\widetilde g_s(\widetilde y^{n}_s,\widetilde z^{n}_s)ds-\int_t^T\widetilde z^{n}_sdW_s+\widetilde k_T^{n,-}-\widetilde k_t^{n,-}-\widetilde k_T^{n,+}+\widetilde k_t^{n,+}.$$

Since we know by Lepeltier and Xu \cite{lx2}, that the above penalization procedure converges to a solution of the corresponding RBSDE, existence and uniqueness are then simple generalization.
\ep

\vspace{0.35em}
We also have a comparison theorem in this context

\begin{Proposition}\label{prop.compdoubly}
Let $\xi_1$ and $\xi_2\in L^2(\mathbb P)$, $V^i$, $i=1,2$ be two adapted, c\`adl\`ag processes and $g^i_s(\omega,y,z)$ two functions, which verify the above assumptions. Let $(y^i,z^i,k^{i,+},k^{i,-})\in \mathbb D^{2}(\mathbb P)\times\mathbb H^{2}(\mathbb P)\times\mathbb I^2(\mathbb P)\times\mathbb I^2(\mathbb P)$, $i=1,2$ be the solutions of the following DRBSDEs with upper obstacle $S^i$ and lower obstacle $L^i$
$$y_t^i=\xi^i+\int_t^Tg^i_s(y^i_s,z^i_s)ds-\int_t^Tz^i_sdW_s+k_T^{i,-}-k_t^{i,-}-k_T^{i,+}+k_t^{i,+}+V_T^i-V_t^i,\text{ }\mathbb P-a.s., \text{ }i=1,2,$$
respectively. If it holds $\mathbb P-a.s.$ that $\xi_1\geq \xi_2$, $V^1-V^2$ is non-decreasing, $S^1\leq S^2$, $L^1\geq L^2$ and $g^1_s(y^1_s,z^1_s)\geq g^2_s(y^1_s,z^1_s)$, then we have for all $t\in [0,T]$
$$Y_t^1\geq Y_t^2,\ \mathbb P-a.s.$$

Besides, if $S^1=S^2$ (resp. $L^1=L^2$), then we also have $dk^{1,+}\geq dk^{2,+}$ (resp. $dk^{1,-}\leq dk^{2,-}$).
\end{Proposition}

\proof
The first part is classical, whereas the second one comes from the fact that the penalization procedure converges in this framework, as seen previously. Indeed, with the notations of the proof of Proposition \ref{ex}, we have in the sense of weak limits
$$k_t^{i,+}=\underset{n\rightarrow +\infty}{\lim}\ n\int_0^t\left(S_s-y_s^{n,i}\right)^-ds.$$

Moreover, using the classical comparison Theorem for RBSDEs with the same lower obstacle, we know that $y^{n,1}\geq y^{n,2}$ and $dk_t^{n,1,-}\leq dk_t^{n,2,-}$. This implies that $dk_t^{n,1,+}\geq dk_t^{n,2,+}$. Passing to the limit yields the result.
\ep

\vspace{0.35em}
Of course, all the above still holds if $\tau$ is replaced by some bounded stopping time $\tau.$ Following Peng's original ideas, we now define a notion of doubly reflected $g$-(super)solutions.
\begin{Definition}
If $y$ is a solution of a DRBSDE of the form \reff{drbsde2}, then we call $y$ a doubly reflected $g$-supersolution on $[0,\tau]$. If $V=0$ on $[0,\tau]$, then we call $y$ a doubly reflected $g$-solution.
\end{Definition}

We have the following proposition concerning the uniqueness of a decomposition of the form \reff{drbsde2}. Notice that unlike in the lower reflected case considered in \cite{mpz}, the processes $V$, $k^+$ and $k^-$ are not necessarily unique.
\begin{Proposition}\label{uniquenessd}
Given $y$ a $g$-supersolution on $[0,\tau]$, there is a unique $z\in\mathbb H^2(\mathbb P)$ and a unique couple $(k^+,k^-,V)\in(\mathbb I^2(\mathbb P))^3$ (in the sense that $V-k^++k^-$ is unique), such that $(y,z,k^+,k^-,V)$ satisfy \reff{drbsde2}.
\end{Proposition}
\proof
If both $(y,z,k^+,k^-,V)$ and $(y,z^1,k^{+,1},k^{-,1},V^1)$ satisfy \reff{drbsde2}, then applying It\^o's formula to $(y_t-y_t)^2$ gives immediately that $z=z^1$ and thus $V-k^++k^-=V^1-k^{+,1}+k^{-,1}$, $\mathbb P-a.s.$
%
\ep

%
%
%

\begin{Remark}
We emphasize once more that the situation here is fundamentally different from \cite{mpz}, where reflected $g$-supersolution were defined for lower reflected BSDEs. In our case, instead of having to deal with the sum of two non-decreasing processes, we actually have to add another non-increasing process. This will raise some difficulties later on, notably when we will prove a non-linear Doob-Meyer decomposition.
\end{Remark}

\subsection{Doob-Meyer decomposition}
We now introduce the notion of doubly reflected $g$-(super)martingales.

\begin{Definition}\label{defdoob}
\begin{itemize}
\item[$\rm{(i)}$] A doubly reflected $g$-martingale on $[0,T]$ is a doubly reflected $g$-solution on $[0,T]$.
\item[$\rm{(ii)}$] A process $(Y_t)$ such that $Y_t\leq S_t$ is a doubly reflected $g$-supermartingale in the strong (resp. weak) sense if for all stopping time $\tau\leq T$ (resp. all $t\leq T$), we have $\mathbb E^\mathbb P[\abs{Y_\tau}^2]<+\infty$ (resp. $\mathbb E^\mathbb P[\abs{Y_t}^2]<+\infty$) and if the doubly reflected $g$-solution $(y_s)$ on $[0,\tau]$ (resp. $[0,t]$) with terminal condition $Y_\tau$ (resp. $Y_t$) verifies $y_\sigma\leq Y_\sigma$ for every stopping time $\sigma\leq\tau$ (resp. $y_s\leq Y_s$ for every $s\leq t$).
\end{itemize}
\end{Definition}

\begin{Remark}
The above Definition differs once more from the one given in \cite{mpz}. Indeed, when defining reflected $g$-supermartingale with a lower obstacle, there is no need to precise that $Y$ is above the barrier $L$, since it is implied by definition (since $y$ is already above the barrier). However, with an upper obstacle, this is not the case and this needs to be a part of the definition.
\end{Remark}

As usual, under mild conditions, a doubly reflected $g$-supermartingale in the weak sense corresponds to a doubly reflected $g$-supermartingale in the strong sense. Besides, thanks to the comparison Theorem, it is clear that a doubly reflected $g$-supersolution on $[0,T]$ is also a doubly reflected $g$-supermartingale in the weak and strong sense on $[0,T]$. The following Theorem addresses the converse property, which gives us a non-linear Doob-Meyer decomposition.

\begin{Theorem}\label{doobmeyerd}
Let $(Y_t)$ be a right-continuous doubly reflected $g$-supermartingale on $[0,T]$ in the strong sense with
$$\mathbb E^\mathbb P\left[\underset{0\leq t\leq T}{\sup}\abs{Y_t}^2\right]<+\infty.$$

Then $(Y_t)$ is a doubly reflected $g$-supersolution on $[0,T]$, that is to say that there exists a quadruple $(z,k^+,k^-,V)\in\mathbb H^2(\mathbb P)\times\mathbb I^2(\mathbb P)\times\mathbb I^2(\mathbb P)\times\mathbb I^2(\mathbb P)$ such that
\begin{align}
\label{rbsde6d}
\begin{cases}
\displaystyle Y_t=Y_T+\int_{t}^T g_s(Y_s,z_s)ds+V_T-V_{t}+k_T^--k_t^--k_T^++k_{t}^+-\int_{t}^T z_sdW_s\\[0.5em]
\displaystyle L_t\leq Y_t\leq S_t, \text{ }\mathbb P-a.s.\\[0.5em]
\displaystyle \int_0^T\left(S_{s^-}-Y_{s^-}\right)dk^+_s=\int_0^T\left(Y_{s^-}-L_{s^-}\right)dk_s^-=0.
\end{cases}
\end{align}
\end{Theorem}

\vspace{0.35em}
We then have the following easy generalization of Theorem $3.1$ of \cite{px} which will be crucial for our proof.

\begin{Theorem}\label{th.monlimd}
Consider a sequence of doubly reflected $g$-supersolutions
\begin{equation}
\label{rbsde5d}
\begin{cases}\displaystyle \tilde y_t^n=\xi+\int_{t}^T g_s(\tilde y_s^n,\tilde z_s^n)ds+\tilde V_T^n-\tilde V_{t}^n+\tilde k_T^{n,-}-\tilde k_t^{n,-}-\tilde k_T^{n,+}+\tilde k_{t}^{n,+}-\int_{t}^T \tilde z_s^ndW_s\\[0.5em]
\displaystyle  L_t\leq\tilde y_t^n\leq S_t\\[0.5em]
\displaystyle \int_0^T\left(S_{s^-}-\tilde y_{s^-}^n\right)d\tilde k_s^{n,+}=\int_0^T\left(\tilde y_{s^-}^n-L_{s^-}\right)d\tilde k_s^{n,-}=0,
\end{cases}
\end{equation}
where the $\tilde V^n$ are in addition supposed to be continuous. Assume furthermore that
 \begin{itemize}
    \item $(\tilde y^n)_{n\geq 0}$ increasingly converges to $\tilde y$ with $\mathbb E^\mathbb P\left[\underset{0\leq t\leq T}{\sup}\abs{\tilde y_t}^2\right]<+\infty.$
    \item $d\tilde k_t^{n,+}\leq d\tilde k_t^{p,+}$ for $n\leq p$ and $\tilde k^{n,+}$ converges to $\tilde k^+$ with $\mathbb E^\mathbb P\left[\left(\tilde k_T^+\right)^2\right]<+\infty.$
    \item $d\tilde k_t^{n,-}\geq d\tilde k_t^{p,-}$ for $n\leq p$ and $\tilde k^{n,-}$ converges to some $\tilde k^-$.
    \item $(\tilde z^n)_{n\geq 0}$ weakly converges in $\mathbb H^2(\mathbb P)$ (along a subsequence if necessary) to $\tilde z$.
\end{itemize}

Then $\tilde y$ is a doubly reflected $g$-supersolution, that is to say that there exists $\tilde V\in \mathbb I^2(\mathbb P)$ such that
$$\begin{cases}
\nonumber\tilde y_t=\xi+\int_{t}^T g_s(\tilde y_s,\tilde z_s)ds+\tilde V_T-\tilde V_{t}+k_T^--k_t^--\tilde k_T^++\tilde k_{t}^+-\int_{t}^T \tilde z_sdW_s,\text{ } \mathbb P-a.s.\\
\nonumber\tilde y_t\leq S_t, \text{ }\mathbb P-a.s.\\
\int_0^T\left(S_{s^-}-\tilde y_{s^-}\right)d\tilde k_s=0, \text{ }\mathbb P-a.s., \text{ }\forall t\in[0,T],
\end{cases}$$

Besides, $\tilde z$ is strong limit of $\tilde z^n$ in $\mathbb H^p(\mathbb P)$ for $p<2$ and $\tilde V_t$ is the weak limit of $\tilde V^n_t$ in $L^2(\mathbb P)$.
\end{Theorem}

\proof
All the convergences are proved exactly as in Theorem $3.1$ of \cite{px}, using the fact that the sequence of added increasing process $\tilde k^{n,-}$ is decreasing. Moreover, since the sequence $y^n$ is increasing it is clear that we have
$$L_t\leq\tilde y_t\leq S_t,\ t\in[0,T],\ \mathbb P-a.s.$$

We now want to show that we also recover the Skorohod conditions. The proofs being similar, we will only show one of them. We have
\begin{align*}
0\leq \int_0^T\left(S_{s^-}-\tilde y_{s^-}\right)d\tilde k_s&=\int_0^T\left(S_{s^-}-\tilde y_{s^-}^n\right)d\tilde k_s+\int_0^T\left(\tilde y^n_{s^-}-\tilde y_{s^-}\right)d\tilde k_s\\
&=\int_0^T\left(S_{s^-}-\tilde y_{s^-}^n\right)d\left(\tilde k_s-\tilde k_s^n\right)+\int_0^T\left(\tilde y^n_{s^-}-\tilde y_{s^-}\right)d\tilde k_s\\
&\leq \int_0^T\left(S_{s^-}-\tilde y_{s^-}^0\right)d\left(\tilde k_s-\tilde k_s^n\right)+\int_0^T\left(\tilde y^n_{s^-}-\tilde y_{s^-}\right)d\tilde k_s
\end{align*}

By the convergences assumed on $\tilde y^n$ and $\tilde k^n$, the right-hand side above clearly goes to $0$ as $n$ goes to $+\infty$, which gives us the desired result.
\ep

\vspace{0.35em}
Let now $Y$ be a given doubly reflected $g$-supermartingale. We follow again \cite{pengg} and we will apply the above Theorem to the following sequence of DRBSDEs
\begin{align}
\label{rbsde7d}
\begin{cases}
\displaystyle y_t^n=Y_T+\int_{t}^T \left(g_s(y_s^n,z_s^n)+n(y_s^n-Y_s)^-\right)ds+k_T^{n,-}-k_t^{n,-}-k_T^{n,+}+k_{t}^{n,+}-\int_{t}^T z_s^ndW_s\\[0.5em]
\displaystyle L_t\leq y_t^n\leq S_t\\[0.5em]
\displaystyle \int_0^T\left(S_{s^-}-y_{s^-}^n\right)dk_s^{n,+}=\int_0^T\left(y_{s^-}^n-L_{s^-}\right)dk_s^{n,-}=0.
\end{cases}
\end{align}

Our first result is that
\begin{Lemma}
For all $n$, we have $Y_t\geq y_t^n.$
\end{Lemma}

\proof
The proof is exactly the same as the proof of Lemma $3.4$ in \cite{pengg}, so we omit it.
\ep

\vspace{0.35em}
We will now prove some estimates which will allow us to apply Theorem \ref{th.monlimd}.

\begin{Lemma}\label{lemd}
There exists a constant $C>0$ independent of $n$ such that the processes defined in \reff{rbsde7d} verify
$$\mathbb E^\mathbb P\left[\underset{0\leq t\leq T}{\sup}\abs{y_t^n}^2+\int_0^T\abs{z_s^n}^2ds+(V_T^n)^2+(k_T^{n,+})^2+(k_T^{n,-})^2\right]\leq C.$$
\end{Lemma}

\proof
First of all, let us define $(\bar y,\bar z,\bar k^+,\bar k^-)$ the unique solution of the DRBSDE with terminal condition $Y_T$, generator $g$, upper obstacle $S$ and lower obstacle $L$ (once again, existence and uniqueness are ensured by the results of \cite{crep} or \cite{lx2}). By the comparison Theorem \ref{prop.compdoubly}, it is clear that we have for all $n\geq 0$
$$y_t^n\geq \bar y_t,\ t\in[0,T],\ \mathbb P-a.s.$$

Consider now $(\tilde y,\tilde z,\tilde k^+,\tilde k^-)$ the unique solution of the doubly reflected BSDE with terminal condition $Y_T$, generator $g$, upper obstacle $S$ and lower obstacle $Y$, that is to say
\begin{align}
\label{drbsde77}
\begin{cases}
\displaystyle \tilde y_t=Y_T+\int_{t}^T g_s(\tilde y_s,\tilde z_s)ds+\tilde k_T^--\tilde k_t^--\tilde k_T^++\tilde k_{t}^+-\int_{t}^T \tilde z_sdW_s\\[0.5em]
\displaystyle Y_t\leq \tilde y_t\leq S_t\\[0.5em]
\displaystyle \int_0^T\left(S_{s^-}-\tilde y_{s^-}\right)d\tilde k_s^+=\int_0^T\left(\tilde y_{s^-}-Y_{s^-}\right)d\tilde k_s^-=0.
\end{cases}
\end{align}

Notice that since the upper obstacle $S$ is a semimartingale  satisfying $$\mathbb E^\mathbb P\left[\underset{0\leq t\leq T}{\sup}\left(\left(S_t\right)^-\right)^2\right]<+\infty,$$ we know from the results of Cr\'epey and Matoussi \cite{crep} (see Theorem $3.2$ and Proposition $5.2$) that the above doubly reflected BSDE has indeed a unique solution and that we have for some constant $C>0$
$$\mathbb E^\mathbb P\left[\left(\tilde k_T^+\right)^2\right]\leq C.$$

Moreover, it is clear that since $\tilde y_s\geq Y_s$, we also have
\begin{equation}
\tilde y_t=Y_T+\int_{t}^T g_s(\tilde y_s,\tilde z_s)ds+n\int_t^T(\tilde y_s-Y_s)^-ds+\tilde k_T^--\tilde k_t^--\tilde k_T^++\tilde k_{t}^+-\int_{t}^T \tilde z_sdW_s.
\label{eq:drbsded}
\end{equation}

Notice also that since $Y$ is a doubly reflected $g$-supermartingale, we have $y\geq L$. Thus we can use now the comparison Theorem of Proposition \ref{prop.compdoubly} for doubly reflected BSDEs with the same upper obstacles. We deduce that
$$y_t^n\leq \tilde y_t,\text{ and }dk_t^{n,+}\leq d\tilde k_t^+.$$

Hence, this implies immediately that for some constant $C$ independent of $n$
\begin{equation}
\mathbb E^\mathbb P\left[\underset{0\leq t\leq T}{\sup}\abs{y_t^n}^2+\left(k_T^{n,+}\right)^2\right]\leq \mathbb E^\mathbb P\left[\underset{0\leq t\leq T}{\sup}\abs{\bar y_t}^2+\underset{0\leq t\leq T}{\sup}\abs{\tilde y_t}^2+\left(\tilde k_T^+\right)^2\right]\leq C.
\label{eq:estimd}
\end{equation}

Define then $V^n_t:=n\int_0^t\left(y_s^n-Y_s\right)^-ds$. We have
\begin{align}\label{hhhd}
\nonumber V_T^n+k_T^{n,-}&=y_0^n-y_T^n-\int_0^Tg_s(y_s^n,z_s^n)ds+k_T^{n,+}+\int_0^Tz_s^ndW_s\\
&\leq C\left(\underset{0\leq t\leq T}{\sup}\abs{y_t^n}+\int_0^T\abs{z_s^n}ds+\int_0^T\abs{g_s(0,0)}ds+k_T^n +\abs{\int_0^Tz_s^ndW_s}\right).
\end{align}

Using \reff{eq:estimd} and BDG inequality, we obtain from \reff{hhhd}
\begin{align}\label{hhhhhd}
\mathbb E^\mathbb P\left[(V_T^n)^2+(k_T^{n,-})^2\right]&\leq C_0\left(1+\mathbb E^\mathbb P\left[\int_0^T\abs{g_s(0,0)}^2ds+\int_0^T\abs{z_s^n}^2ds\right]\right).
\end{align}

Then, using It\^o's formula, we obtain classically for all $\epsilon>0$
\begin{align*}
\nonumber\mathbb E^\mathbb P\left[\int_0^T\abs{z_s^n}^2ds\right]&\leq\mathbb E^\mathbb P\left[(y_T^n)^2+2\int_0^Ty_s^ng_s(y_s^n,z_s^n)ds +2\int_0^Ty_{s^-}^nd(V_s^n-k_s^{n,+}+k_s^{n,-})\right]\\
&\leq\mathbb E^\mathbb P\left[C(1+\underset{0\leq t\leq T}{\sup}\abs{y_t^n}^2)+\int_0^T\frac{\abs{z_s^n}^2}{2}ds+\epsilon(|V_T^n|^2+|k_T^{n,+}|^2+|k_T^{n,-}|^2)\right].
\end{align*}

Then, from \reff{hhhhhd}, we obtain by choosing $\epsilon=\frac{1}{4C_0}$ that
$$\mathbb E^\mathbb P\left[\int_0^T\abs{z_s^n}^2ds\right]\leq C.$$

Reporting this in \reff{hhhhhd} ends the proof.
\ep

\vspace{0.95em}
Finally, we can now prove Theorem \ref{doobmeyerd}

\vspace{0.35em}
\proof[Proof of Theorem \ref{doobmeyerd}]
We first notice that since $Y_t\geq y_t^n$ for all $n$, by the comparison Theorem for DRBSDEs, we have $$y_t^n\leq y_t^{n+1}, \ dk_t^{n,-}\geq dk_t^{n+1,-} \text{ and } dk_t^{n,+}\leq dk_t^{n+1,+}.$$ By the a priori estimates of Lemma \ref{lemd}, they therefore converge to some processes $y$, $k^+$ and $k^-$. Moreover, since $z^n$ is bounded uniformly in $n$ in the Banach $\mathbb H^2(\mathbb P)$, there exists a weakly convergent subsequence, and the same holds for $g_t(y_t^n,z_t^n)$. Hence, all the conditions of Theorem \ref{th.monlimd} are satisfied and $y$ is a doubly reflected $g$-supersolution on $[0,T]$ of the form

$$y_t=Y_T+\int_t^Tg_s(y_s,z_s)ds+V_T-V_t-k_T^++k_t^++k_T^--k_t^--\int_t^Tz_sdW_s,$$
where $V_t$ is the weak limit of $V_t^n:=n\int_0^t(y_s^n-Y_s)^-ds$. From Lemma \ref{lemd}, we have
$$\mathbb E^\mathbb P[(V_T^n)^2]=n^2\mathbb E^\mathbb P\left[\int_0^T\abs{(y_s^n-Y_s)^-}^2ds\right]\leq C.$$

It then follows that $Y_t=y_t$ (since we already had $Y_t\geq y_t^n$ for all $n$), which ends the proof.
\ep

\subsection{Time regularity of doubly reflected $g$-supermartingales}

In this section we prove a downcrossing inequality for doubly reflected $g$-supermartingales in the spirit of the one proved in \cite{cpeng}. We use the same notations as in the classical theory of $g$-martingales (see \cite{cpeng} and \cite{pengg} for instance).

\begin{Theorem}\label{down}
Assume that $g(0,0)=0$. Let $(Y_t)$ be a positive reflected $g$-supermartingale in the weak sense and let $0=t_0<t_1<...<t_i=T$ be a subdivision of $[0,T]$. Let $0\leq a<b$, then there exists $C>0$ such that $D_a^b[Y,n]$, the number of downcrossings of $[a,b]$ by $\left\{Y_{t_j}\right\}$, verifies
$$\mathcal E^{-\mu}[D_a^b[Y,n]]\leq\frac{C}{b-a}\mathcal E^\mu[Y_0\wedge b],$$
where $\mu$ is the Lipschitz constant of $g$.
\end{Theorem}

\proof
Consider for all $i=0..n$ the following DRBSDEs with upper obstacle $S$ and lower obstacle $L$ on $[0,t_i]$
$$\begin{cases}
\displaystyle y_t^i=Y_{t_i}-\int_t^T(\mu\abs{y_s^i}+\mu\abs{z_s^i})ds+k_{t_i}^{i,-}-k_{t}^{i,-}-k_{t_i}^{i,+}+k_{t}^{i,+}-\int_{t}^{t_i} z_s^idW_s\text{ } \mathbb P-a.s.\\[0.5em]
\displaystyle L_t\leq y_t^i\leq S_t, \text{ }\mathbb P-a.s.\\[0.5em]
\displaystyle \int_0^{t_i}\left(y_{s^-}^i-L_{s^-}\right)dk_s^{i,-}=\int_0^{t_i}\left(S_{s^-}-y_{s^-}^i\right)dk_s^{i,+}=0, \text{ }\mathbb P-a.s.
\end{cases}$$

By the comparison theorem of Proposition \ref{prop.compdoubly}, we know that we have for all $i$, $y_t^i\geq \tilde y_t^i$ for $t\in[0,t_i]$, where $(\tilde y^i,\tilde z^i,\tilde k^i)$ is the unique solution of the RBSDE on $[0,t_i]$ with the same generator and terminal condition as above and upper obstacle $S$, that is to say
$$\begin{cases}
\tilde y_t^i=Y_{t_i}-\int_t^T(\mu\abs{\tilde y_s^i}+\mu\abs{\tilde z_s^i})ds-\tilde k_{t_i}^{i}+\tilde k_{t}^{i}-\int_{t}^{t_i} \tilde z_s^idW_s\text{ } \mathbb P-a.s.\\[0.5em]
\tilde y_t^i\leq S_t, \text{ }\mathbb P-a.s.\\[0.5em]
\int_0^{t_i}\left(S_{s^-}-\tilde y_{s^-}^i\right)d\tilde k_s^{i}=0, \text{ }\mathbb P-a.s.
\end{cases}$$

We define $a_s^i:=-\mu\ \text{sgn}(z_s^i)1_{t_{j-1}<s\leq t_j}$ and $a_s:=\sum_{i=0}^na_s^i$. Let $\mathbb Q^a$ be the probability measure defined by
$$\frac{d\mathbb Q^a}{d\mathbb P}=\mathcal E\left(\int_0^Ta_sdW_s\right).$$

\vspace{0.35em}
We then have easily that $y_t^i\geq 0$ since $Y_{t_i}\geq 0$. We next define $\hat y^i:=\tilde y^i-k^i,\ \hat Y^i:=Y-k^i.$ Then, $(\hat y^i,z^i)$ solves the following BSDE on $[0,t_i]$
$$\hat y^i_t=\hat Y^i_{t_i}-\int_t^{t_i}\mu\left(\hat y^i_s+k_s^i\right)+\mu\abs{z^i_s}ds-\int_t^{t_i}z_s^idW_s.$$

It is then easy to solve this BSDE to obtain
$$\hat y_t^i=\mathbb E_t^{\mathbb Q^a}\left[e^{-\mu(t_i-t)}\hat Y^i_{t_i}-\mu\int_t^{t_i}e^{-\mu(s-t)}k^i_sds\right].$$

Define now the following c\`adl\`ag process $k_t:=\sum_{i=i}^nk_t^i1_{t_{i-1}\leq t<t_i}$ and $\hat Y:=Y-k$. We clearly have for $t=t_{i-1}$
$$\hat y_{t_{i-1}}^i=\mathbb E_{t_{i-1}}^{\mathbb Q^a}\left[e^{-\mu(t_i-t_{i-1})}\hat Y_{t_i}-\mu\int_{t_{i-1}}^{t_i}e^{-\mu(s-t_{i-1})}k_sds\right].$$

Now, since $Y$ is a doubly reflected $g$-supermartingale (and thus also a doubly reflected $g^{-\mu}$-supermartingale where $g^{-\mu}_s(y,z):=-\mu(\abs{y}+\abs{z})$ by a simple application of the comparison theorem), we have
$$\hat y^i\leq y^i-k^i\leq \hat Y.$$

Hence, we have obtained
$$\mathbb E_{t_{i-1}}^{\mathbb Q^a}\left[e^{-\mu(t_i-t_{i-1})}\hat Y_{t_i}-\mu\int_{t_{i-1}}^{t_i}e^{-\mu(s-t_{i-1})}k_sds\right]\leq  \hat Y_{t_{i-1}}.$$

This actually implies that the process $X:=(X_{t_i})_{0\leq i\leq n}$ where
$$X_{t_i}:=e^{-\mu t_{i}}\hat Y_{t_i}-\mu\int_0^{t_i}e^{-\mu s}k_sds,$$
is a $\mathbb Q^a$-supermartingale. Then we can finish the proof exactly as in the proof of Theorem $6$ in \cite{cpeng}.
\ep

\end{appendix}

\end{document}